\documentclass[sn-mathphys,Numbered]{sn-jnl}

\RequirePackage{xcolor}[2.11]
\definecolor{header1}{cmyk}{.9,.5,0,.35}
\definecolor{blue1}{cmyk}{.9,.7,0,0}
\definecolor{blue2}{cmyk}{.93,.95,.2,.07}
\definecolor{maroon}{cmyk}{.4,1,.3,.2}
\definecolor{gold1}{cmyk}{.2,.2,1,.1}
\definecolor{gray}{cmyk}{0,0,0,.5}
\definecolor{green1}{cmyk}{1,0,1,0}
\definecolor{proofcolor}{cmyk}{1,0,1,0}
\definecolor{red1}{cmyk}{0,1,.8,0}
\definecolor{orange1}{cmyk}{0,.55,1,0}
\definecolor{strip}{cmyk}{.6,.1,.1,.1}

\RequirePackage[amsmath,thmmarks,hyperref]{ntheorem}[1.33]

\newcommand{\proofbox}{\textcolor{header1}{\rule{1.5ex}{1.5ex}}}

\theorempreskip{1ex plus .25ex minus .1ex}
\theorempostskip{1ex plus .25ex minus .1ex}

\theoremstyle{nonumberplain}
\theoremheaderfont{\normalfont\itshape}
\theorembodyfont{\normalfont}
\theoremseparator{.}
\theoremsymbol{\proofbox}
\newtheorem{proof}{Proof}

\theoremstyle{plain}
\theoremheaderfont{\normalfont\sffamily}
\theorembodyfont{\normalfont\itshape}
\theoremseparator{.}
\theoremsymbol{}

\newcommand{\newsiamthm}[2]{
  \theoremstyle{plain}
  \theoremheaderfont{\normalfont\sffamily}
  \theorembodyfont{\normalfont\itshape}
  \theoremseparator{.}
  \theoremsymbol{}
  \newtheorem{#1}[theorem]{#2}
}

\newsiamthm{proposition}{Proposition}
\newtheorem{prop}{Proposition}

\usepackage{graphicx}%
\usepackage{multirow}%
\usepackage{amssymb,amsfonts}%
\usepackage{bm}
\usepackage{mathtools}
\usepackage{mathrsfs}%
\usepackage[title]{appendix}%
\usepackage{xcolor}%
\usepackage{textcomp}%
\usepackage{manyfoot}%
\usepackage{booktabs}%
\usepackage{algorithm}%
\usepackage{algorithmicx}%
\usepackage{algpseudocode}%
\usepackage{listings}%
\usepackage{caption}
\usepackage{subcaption}

\usepackage{enumitem}
\newlist{propenum}{enumerate}{1}
\setlist[propenum,1]{%
  label = \upshape(\alph*),
  ref = {\alph*},
}
\crefname{prop}{Proposition}{Propositions}
\crefname{propenumi}{Proposition}{Propositions}
\crefformat{propenumi}{Proposition~(#2#1#3)}
\newcommand\crefproppart[2]{\namecref{#1}~\labelcref{#1}(\ref{#2})}

\begin{document}

\title[Article Title]{Structure of Periodic Orbit Families in the Hill Restricted 4-Body Problem}

\author*[1]{\fnm{Gavin M.} \sur{Brown}}\email{gavin.m.brown@colorado.edu}

\author[1]{\fnm{Luke T.} \sur{Peterson}}\email{Luke.Peterson@colorado.edu}

\author[1]{\fnm{Damennick B.} \sur{Henry}}\email{Damennick.Henry@colorado.edu}

\author[1]{\fnm{Daniel J.} \sur{Scheeres}}\email{scheeres@colorado.edu}

\affil[1]{\orgdiv{Ann and H.J. Smead Aerospace Engineering Sciences}, \orgname{University of Colorado Boulder}, \state{CO}, \country{USA}}

\abstract{The Hill Restricted 4-Body Problem~(HR4BP) is a coherent time-periodic model that can be used to represent motion in the Sun-Earth-Moon~(SEM) system. Periodic orbits were computed in this model to better understand the periodic orbit family structures that exist in these types of systems. First, periodic orbits in the Circular Restricted 3-Body Problem~(CR3BP) representation of the Earth-Moon~(EM) system were identified. A Melnikov-type function was used to identify a set of candidate points on the EM CR3BP periodic orbits to start a continuation algorithm. A pseudo-arclength continuation scheme was then used to obtain the corresponding periodic orbit families in the HR4BP when including the effect of the Sun. Bifurcation points were identified in the computed families to obtain additional orbit families.}

\keywords{Periodic orbits, Bifurcations, Three-Body Problems, Melnikov Function, Sun-Earth-Moon System}

\maketitle

\newpage 

\tableofcontents

\newpage

\section{Introduction} \label{s:INTR}
Cislunar space is the operating environment for many current and future spacecraft missions. It is the target of NASA's planned Gateway, a crewed space station in orbit near the Moon, and the focus of many other NASA objectives~\cite{cit:NASASP22,cit:WhitleyDavisBurkeEtal}. The Circular Restricted 3-Body Problem~(CR3BP) is one of the most basic dynamical models we have to represent motion in cisulunar space.
There are several dynamical models that seek to present a more realistic representation of motion in cislunar space than the CR3BP. The Elliptic Restricted 3-Body Problem~(ER3BP) includes the effect of the eccentricity of the two primaries' orbit~\cite{cit:Szebehely,cit:ParkHowell,cit:PengBai}. Other models also attempt to account for the effect of the Sun on the motion of a spacecraft in the Earth-Moon~(EM) system. One such model is the Bicircular Restricted 4-Body Problem~(BCP), which is incoherent as it does not account for the effect of the Sun on the Earth or the Moon~\cite{cit:Huang,cit:GomezJorbaMasdemontEtal,cit:Rosales}. In addition, there is no accurate dynamical equivalent to $L_2$ which is one major drawback of the BCP~\cite{cit:JorbaCuscoFarresJorba}. Periodic orbits in the BCP corresponding to the CR3BP triangular points have been computed previously~\cite{cit:SimoGomezJorbaEtal}, as have families of periodic orbits~\cite{cit:BoudadHowellDavis} and quasi-periodic orbits~\cite{cit:CastellaJorba,cit:RosalesJorbaJorbaCuscoBCP}. Unlike the BCP, the Quasi-Bicircular Model~(QBCP) accounts for the effect of the Sun on the Earth and Moon by modeling their motion as a solution to the 3-body problem~\cite{cit:Andreu}. Previously, quasi-periodic orbits~(QPOs) have been computed in this coherent model~\cite{cit:Andreu,cit:RosalesJorbaJorbaCuscoQBCP}.

The Hill Restricted 4-Body Problem~(HR4BP), developed by Scheeres~\cite{cit:ScheeresHR4BP} in 1998, is another coherent time periodic model describing the motion of a small body~($P_3$) in the presence of three large bodies~($P_0$, $P_1$, and $P_2$). The model is a higher fidelity model than the CR3BP and BCP, more accurately represents the true dynamics in the Sun-Earth-Moon~(SEM) system, is easier to implement than the QBCP, and has previously been used to study the Sun's effect on the EM system~\cite{cit:ScheeresHR4BP,cit:PetersonRosalesScheeres,cit:OlikaraGomezMasdemont}. Specifically, periodic orbits and quasi-periodic orbits~(QPOs) related to the EM libration points have been computed~\cite{cit:OlikaraGomezMasdemont,cit:HenryRosalesBrownEtalISTS,cit:HenryRosalesBrownEtalAASC,cit:PetersonJorbaBrownEtal}, as have connections between Sun-Earth and EM libration point orbits~\cite{cit:OlikaraScheeres}. Olikara, G{\'o}mez, and Masdemont \cite{cit:OlikaraGomezMasdemont} computed the HR4BP periodic orbit families corresponding to the CR3BP EM $L_1$ and $L_2$ points, the $L_1$ Lyapunov orbit with $T = \pi$, as well as two additional periodic orbit families existing near a broken pitchfork bifurcation in the dynamical equivalent of the EM $L_2$ point. Sanaga and Howell~\cite{cit:SanagaHowell} computed several orbit families in the HR4BP corresponding to synodic resonant halo orbits. In that work, the authors presented an additional dynamical model, a reduced model of the HR4BP referred to as the RHR4BP~\cite{cit:SanagaHowell}. In this work, by using a methodology that leverages a Melnikov-type function and techniques related to detecting bifurcation points, we will present a more detailed picture of the periodic orbits that exist in the HR4BP.

\section{Problem Formulation} \label{s:PRFO}
\subsection{Formulation of the HR4BP}
\label{sub:H4EM}
\begin{figure}[!ht]
    \centering
    \includegraphics[width=0.6\linewidth]{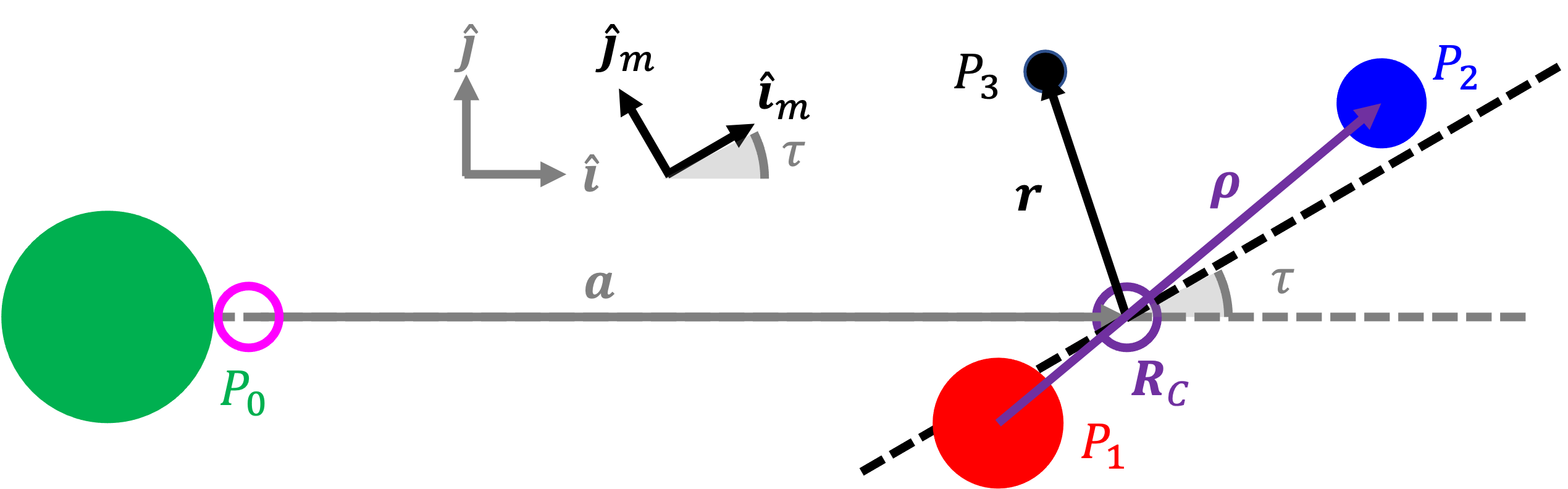}
    \caption{Coordinate Frames Used in~{\rm{\cite{cit:ScheeresHR4BP}}}. Note $\bm{\rho} = \hat{\bm{\imath}}_m + \bar{\bm{\rho}}$.}
    \label{fig:HR4BPCF}
\end{figure}

The equations of motion for the HR4BP were originally presented by Scheeres~\cite{cit:ScheeresHR4BP}. Let $P_0$ represent the largest body, $P_1$ and $P_2$ represent the two other bodies with non-negligible mass, and $P_3$ be a body with negligible mass near $P_1$ and $P_2$. $\mathcal{B}: \left\{ \hat{\bm{\imath}}_m, \hat{\bm{\jmath}}_m, \hat{\bm{k}} \right\}$ is defined as a rotating frame with a constant angular velocity (with respect to inertial space) in the direction $\hat{\bm{\Omega}} = \hat{\bm{k}}$ and an origin at the center of mass of $P_1$, $P_2$, and $P_3$. Note $P_0$, $P_1$, and $P_2$ lie in the $xy$-plane (i.e., the plane perpendicular to $\hat{\bm{k}}$). All vector components in this work are presented in the $\mathcal{B}$-frame. Visual representations of this frame, an intermediate frame used in~\cite{cit:ScheeresHR4BP}, and terms defined later in this section are provided in \cref{fig:HR4BPCF}, where the magenta circle is the total system center of mass and the purple circle is the center of mass of $P_1$, $P_2$, and $P_3$.

Let $M_i$ represent the mass of body $P_i$. $P_1$ and $P_2$ are referred to as the primaries. $n_0$ represents the mean motion of $P_0$ considering a 2-body system of $P_0$ as one body and treating the two primaries as one combined body at $\bm{R}_c$. $\mu$ represents the mass ratio of $P_1$ and $P_2$, and $m$ represents the relationship between the period of $\bm{R}_c$ about the total system center of mass and the period of $P_1$ and $P_2$ about $\bm{R}_c$~\cite{cit:OlikaraGomezMasdemont}. $n$ represents the mean motion of the primaries in the 2-body system of $P_1$ and $P_2$. We will assume $\bm{R}_c \approx \bm{a} = d_a \hat{\bm{i}}$ where $d_a$ is the distance between the the Sun and the Earth-Moon barycenter~\cite{cit:ScheeresHR4BP}.

Expressions for $\mu$ and $m$, as well as a third parameter $\nu$, are provided in \cref{eqn:HR4BPPar}~\cite{cit:OlikaraGomezMasdemont,cit:SanagaHowell}. To be consistent with the standard CR3BP notation, $\mu$ and $\nu$ are defined differently here compared to how they are defined in~\cite{cit:ScheeresHR4BP}.
\begin{subequations}
    \begin{align}
        \mu &= \frac{M_2}{M_1 + M_2}
        \label{eqn:HR4BPPamu} \\
        m &= \frac{n_0/n}{1 - n_0/n}
        \label{eqn:HR4BPParnm} \\
        \nu &= \frac{M_1 + M_2}{M_0}
        \label{eqn:HR4BPParnu}
    \end{align}
    \label{eqn:HR4BPPar}
\end{subequations}

In the SEM system, $m$ is the ratio of the difference between the synodic and sidereal months relative to a sidereal month (i.e., $m \approx 0.0808$), and $\mu \approx 0.0122$. Time is scaled such that the relative physical configuration of $P_0$, $P_1$, and $P_2$ repeats every $2 \pi$ time units (i.e., $2 \pi$ time units is equal to a synodic month)~\cite{cit:OlikaraGomezMasdemont}. However, the equations of motion are $\pi$-periodic~\cite{cit:ScheeresHR4BP}. Please refer to \cref{app:HEOM} for more detail on the terms used in the equations of motion. Distance is scaled with $1 \, \text{DU}$ set to be the average distance between $P_1$ and $P_2$. Note that $d_a$ can be computed based on $1 \, \text{DU}$, $\nu$, and $a_0$ (a term related to the particular solution for the motion of the primaries defined in \cref{eqn:HVOa0} in \cref{app:HEOM}) which is a function of $m$. All results will be presented in these scaled units. The scaled mass, distance, and time units are provided in \cref{eqn:HR4BPND}~\cite{cit:ScheeresHR4BP,cit:SanagaHowell}.
\begin{subequations}
    \begin{align}
        1 \; \text{MU} &= M_0
        \label{eqn:HR4BPNDm} \\
        1 \; \text{DU} &= a_0 d_a \nu^{1/3}
        \label{eqn:HR4BPNDd} \\
        1 \; \text{TU} &= m / n_0 = \left( 1 + m \right) / n
        \label{eqn:HR4BPNDt}
    \end{align}
    \label{eqn:HR4BPND}
\end{subequations}

The structure of the HR4BP equations of motion will now be briefly presented where $\bm{r}$ represents the position of $P_3$.
In this work, $'$ indicates a time derivative with respect to scaled time $\tau$.
The positions of the two primaries $P_1$ and $P_2$ (e.g., the Earth and the Moon) when accounting for the effect of a much more massive body $P_0$ (e.g., the Sun) can be obtained by determining the Hill variation orbit coefficients. These coefficients are used in a Fourier series representation of the variation orbit family solution to the Hill equations of motion~\cite{cit:ScheeresHR4BP}. $P_3$'s relative position with respect to $P_1$ and $P_2$ are $\bm{R}_{1 - \mu}$ and $\bm{R}_{\mu}$, respectively.
\begin{subequations}
    \begin{align}
        \bm{R}_{1 - \mu} &= \bm{r} + \mu \left( \hat{\bm{\imath}}_m + \bar{\bm{\rho}} \right)
        \; \text{and} \; R_{1 - \mu} = \left| \bm{R}_{1 - \mu} \right|
        \label{eqn:EOMS1R1mmu} \\
        \bm{R}_{\mu} &= \bm{r} - (1 - \mu) \left( \hat{\bm{\imath}}_m + \bar{\bm{\rho}} \right)
        \; \text{and} \; R_{\mu} = \left| \bm{R}_{\mu} \right|
        \label{eqn:EOMS1Rmu}
    \end{align}
    \label{eqn:EOMS1}
\end{subequations}

The potential ($V$) can be evaluated using the expression in \cref{eqn:HR4BPEOMV}~\cite{cit:ScheeresHR4BP}.
\begin{equation}
    \begin{split}
        V =& \frac{1}{2} \left( 1 + 2 m + \frac{3}{2} m^2\right) \left( x^2 + y^2 \right) - \frac{1}{2} m^2 z^2 + \frac{3}{4} m^2 \left( \left( x^2 - y^2 \right) \cos{2 \tau} - 2 x y \sin{2 \tau} \right) \\
        &+ \frac{m^2}{a_0^3} \left[ \frac{1 - \mu}{R_{1 - \mu}} + \frac{\mu}{R_\mu} \right]
    \end{split}
    \label{eqn:HR4BPEOMV}
\end{equation}

Let $\left[ \Phi \left( \tau , \tau_0 \right) \right] = \frac{\partial \bm{X}}{\partial \bm{X}_0}$ be the state transition matrix, $\left[ \Psi \left( \tau , \tau_0 \right) \right] = \left[ \Psi_m \; \Psi_\mu \right] = \left[ \frac{\partial \bm{X}}{\partial m} \; \frac{\partial \bm{X}}{\partial \mu} \right]$, $\left[ A \right] = \frac{\partial \bm{X}'}{\partial \bm{X}}$ be the Jacobian matrix, and $\left[ C \right] = \left[ \frac{\partial \bm{X}'}{\partial m} \; \frac{\partial \bm{X}'}{\partial \mu} \right]$. Note $\bm{X}_0$ refers to the state at initial time $\tau_0$, $\left[ \Phi \left( \tau_0 , \tau_0 \right) \right] = \left[ I_{6 \times 6} \right]$, and $\left[ \Psi \left( \tau_0 , \tau_0 \right) \right] = \left[ 0_{6 \times 2} \right]$. Under the assumption $\bm{\Omega} = \hat{\bm{k}}$, the equations to integrate $\bm{X}$, $\left[ \Phi \right]$, and $\left[ \Psi \right]$ are:
\begin{subequations}
    \begin{align}
        \bm{X}' &= \begin{bmatrix} \bm{r}' \\ - 2 \left( 1 + m \right) \bm{\Omega} \times \bm{r}' + \nabla V \end{bmatrix}
        \label{eqn:EOMDX} \\
        \left[ \Phi \right]' &= \left[ A \right] \left[ \Phi \right]
        \label{eqn:EOMDPhi} \\
        \left[ \Psi \right]' &= \left[ A \right] \left[ \Psi \right] + \left[ C \right]
        \label{eqn:EOMDPsi}
    \end{align}
    \label{eqn:HR4BPEOM}
\end{subequations}

It is important to note that when $m \rightarrow 0$, the HR4BP equations of motion take on the form of the CR3BP equations of motion. As $m$ increases from zero, the effect of the more massive body on the system becomes more pronounced, and $m_{\text{SEM}} = 0.0808$ for the SEM system.
In addition to being a periodically forced system, there are also important symmetries in the equations of motion~\cite{cit:SanagaHowell,cit:ScheeresHR4BP}.
\begin{subequations}
    \begin{align}
        S_1 &: \left( x, y, z, x', y', z', k \pi + \tau \right) \rightarrow \left( x, -y, z, -x', y', -z', k \pi - \tau \right)
        \label{eqn:HR4BPSym1} \\
        S_2 &: \left( x, y, z, x', y', z', k \pi + \frac{\pi}{2} + \tau \right) \rightarrow \left( x, -y, z, -x', y', -z', k \pi + \frac{\pi}{2} - \tau \right)
        \label{eqn:HR4BPSym2}
    \end{align}
    \label{eqn:HR4BPSym}
\end{subequations}

As an example, from~\cite{cit:ScheeresHR4BP}, if the initial condition $\bm{X}(\tau = \tau_0) = \left[ x_0, y_0, 0, x'_0, y'_0, 0 \right]^T$ corresponds to a periodic orbit with a period that is an integer multiple of $\pi$, then $\bm{X}(\tau = -\tau_0) = \left[ x_0, -y_0, 0, -x'_0, y'_0, 0 \right]^T$ also corresponds to a periodic orbit with the same period. This symmetry will be an important consideration when computing periodic orbits.

\subsection{Melnikov Theory}
\label{sub:MELT}
When a time-periodic perturbation is added to an autonomous system, a single structure in the autonomous system may have multiple ``dynamical equivalents'' in the perturbed system (i.e., related structures that exist in the perturbed system). For example, at least four dynamical equivalents to the 9:2 NRHO in the BCP have been identified previously~\cite{cit:BoudadHowellDavis}. In that system the continuation of these orbits followed many possible paths, and we expect a similarly complicated behavior in the HR4BP. Melnikov theory is one of the tools we can use to analyze these perturbed systems.
For an introduction to classical Melnikov theory, which was originally presented by Melnikov~\cite{cit:Melnikov}, we refer the reader to sections 3.1 and 3.2 of~\cite{cit:GreenspanHolmes}, but additional detail can be found in Chapter~28.1 of~\cite{cit:Wiggins}, Chapter~4.5 of~\cite{cit:GuckenheimerHolmes}, and Chapter~4.9 of~\cite{cit:Perko}.

This theory has been improved, extended, and applied to a variety of different systems, including those with higher dimensions~\cite{cit:Wiggins,cit:Haller,cit:GuoTianXueEtal,cit:VeermanHolmes,cit:Yagasaki,cit:PolcarSemerak}. By considering an expansion of the ``energy principle'' (the change in energy at the start and end of any periodic orbit must be zero), Cenedese and Haller~\cite{cit:CenedeseHaller} present an expression for a Melnikov-type function that is applicable to our study of orbits in the HR4BP. It is this Melnikov-type function of~\cite{cit:CenedeseHaller} that we will use in our analysis, which we will simply refer to as the Melnikov function.
We will now present the basic form of the Melnikov function from~\cite{cit:CenedeseHaller}. Say the acceleration of an autonomous system (in our case the CR3BP) is given by $\ddot{\bm{r}} = \bm{f} \left( \bm{X} \right)$. Let the flow of the vector field in this autonomous system be defined as $\bm{X} \left( \tau \right) = \varphi \left( \bm{X}_0, \tau \right)$. Say a small perturbation is added so that the new acceleration is
\begin{equation}
    \ddot{\bm{r}} = \bm{F} = \bm{f} \left( \bm{X} \right) + \varepsilon \bm{g} \left( \bm{X}, \tau_0 + \tau; T_g, \varepsilon \right)
    \label{eqn:MTEOMmain}
\end{equation}
where $\bm{g}$ is time-periodic with period $T_g$ so $\bm{g} \left( \bm{X}, \tau_0 + \tau; T_g, \varepsilon \right) = \bm{g} \left( \bm{X}, \tau_0 + \tau + T_g; T_g, \varepsilon \right)$, $\varepsilon$ is the perturbation parameter, and $\tau$ corresponds to time relative to $\tau_0$.

Periodic orbits in this system must have a minimal period~($T$) that is in some resonance with the forcing period~($T_g$) where $T_g = \pi$ for the HR4BP (i.e., $T = b T_g$ where $b$ is a positive integer)~\cite{cit:RhoumaChicone}.
Say there is a specific periodic orbit $\mathcal{Z}$ in the unperturbed autonomous system of period $T^*$ which satisfies the appropriate resonance condition with $T_g$ (i.e., $b T_g = a T^*$ where $a$ and $b$ are relatively prime integers), and an initial state time history of the orbit in the unperturbed system is generated using a specific state $\prescript{*}{}{\bm{X}}_s$ on the orbit. As the periodic orbit exists in the autonomous system, any point on that orbit $\prescript{*}{}{\bm{X}} \left( s \right) = \varphi \left( \prescript{*}{}{\bm{X}}_s, s \right)$ where $s \in \left[ 0, T^* \right)$ could be used as the initial state to generate the same periodic orbit, but with a different time history associated with each state. So, based on the particular value of $s$, the periodic orbit state time history in the unperturbed system can be represented as $\prescript{*}{}{\bm{X}} \left( s + \tau \right) = \varphi \left( \prescript{*}{}{\bm{X}}_s, s + \tau \right)$ for $\tau \in \left[ 0, T^* \right]$ where $\prescript{*}{}{\bm{X}} \left( s + \tau \right) = \prescript{*}{}{\bm{X}} \left( s + \tau + T^* \right)$ and $s \in \left[ 0, T^* \right)$. With that in mind, the value of the Melnikov function takes the form:
\begin{equation}
    \mathcal{M} \left( s, \tau_0 \right) = \int_{0}^{a T^*}{\bm{g} \left( \prescript{*}{}{\bm{X}} \left( s + \tau \right), \tau_0 + \tau \right) \cdot \prescript{*}{}{\bm{r}'} \left( s + \tau \right) \, \text{d}\tau}
    \label{eqn:MFmmain}
\end{equation}
which is dependent on the choice of the initial point on the unperturbed periodic orbit (i.e., the value of $s$) and the initial time $\tau_0$. Note that we have removed the $T_g$ and $\varepsilon$ terms from this equation for conciseness, and the period of the orbit in the perturbed system corresponding to the orbit in the unperturbed system is $T = b T_g$.

For an orbit to be a periodic solution, there must be no work done by the perturbation force on the orbit over one period. The zeros of the Melnikov function represent where this is the case (at least when considering the leading order terms). Provided $\mathcal{M} \left( s, \tau_0 \right)$ is not identically zero for all $s \in \left[ 0, T^* \right)$, we expect to be able to continue periodic orbits from the unperturbed system into the perturbed system at the points $s$ on the unperturbed orbit provided the initial time when beginning the integration is $\tau_0$. Points where this function is zero are points where we will attempt to continue CR3BP orbits into the HR4BP.

\section{Methodology} \label{s:METH}
As this system is non-autonomous, we do not expect to identify families of periodic orbits as we would in an autonomous system~\cite{cit:ScheeresBook}. However, we do expect to identify periodic orbit ``families'' in the HR4BP if we fix the value of $T$ and initial time ($\tau_0$), and perform a continuation while allowing the initial state ($\bm{X}_0 = \bm{X}(\tau = \tau_0)$) and at least one of the parameters $m$ and/or $\mu$ to vary.
The phrase ``HR4BP periodic orbit family'' will refer to the ``family'' that can be computed starting from a specific periodic orbit in the HR4BP and then allowing $\bm{X}_0$ and $m$ to vary. To compute these families we will first expand the HR4BP equations of motion about $m = 0$. We will then use this expansion when evaluating the Melnikov function for orbits in the EM CR3BP with appropriate periods. We will use pseudo-arclength continuation starting at $m = 0$ to compute the corresponding HR4BP periodic orbit families. We will start at the five libration points whose HR4BP orbits have periods of $T = \pi$, and at each state on the selected CR3BP periodic orbits where $\mathcal{M}(s,\tau_0) = 0$. We will then identify bifurcations along these families in order to compute additional families starting at non-zero values of $m$. Note we will always use $\tau_0 = 0$ unless explicitly stated otherwise.

\subsection{Expanding the HR4BP Equations of Motion}
\label{sub:EHEM}
In order to use the Melnikov function presented in \cref{eqn:MFmmain}, we need to expand the HR4BP equations of motion in \cref{eqn:HR4BPEOM} about $m = 0$ so they are in the form of \cref{eqn:MTEOMmain}. This expansion will result in $\varepsilon \bm{g} = \bm{g}_1 m + \bm{g}_2 m^2 + \bm{g}_3 m^3 + \cdots$. By using $\varepsilon = m$, the leading order term can be used as $\bm{g}$ in \cref{eqn:MFmmain} (i.e., $\bm{g} = \bm{g}_1$). While the definition of $M$ given in \cref{eqn:HVOM} will be used for all calculations except when evaluating the Melnikov function. In this case, the definition $M = m$ will be used which results in different values for the Hill variation orbit coefficients: $d_0 = 1$, $d_1 = -2/3$, $d_2 = 7/18$, $d_3 = -4/81$, $c_{-1,2} = -19/16$, $c_{1,2} = 3/16$, $c_{-1,3} = -5/3$, and $c_{1,3} = 1/2$. Note that these are the same values that are presented in~\cite{cit:ScheeresHR4BP}. Let $\bm{R}_{1 - \mu,\text{C}}$ and $\bm{R}_{\mu,\text{C}}$ be the relative position of $P_3$ with respect to $P_1$ and $P_2$, respectively, in the CR3BP (i.e., when $m = 0$). Performing this expansion we obtain the following result for $\bm{g}_1$:
\begin{subequations}
    \begin{align}
        \bm{g}_1 &=
        2 \begin{bmatrix}
            y' \\ -x' \\ 0
        \end{bmatrix} +
        2 \begin{bmatrix}
            x \\ y \\ 0
        \end{bmatrix} + 
        3 d_1 \left( \frac{1 - \mu}{R_{1 - \mu,\text{C}}^3} \bm{R}_{1 - \mu,\text{C}} + \frac{\mu}{R_{\mu,\text{C}}^3} \bm{R}_{\mu,\text{C}} \right)
        \label{eqn:HR4BPg1old1} \\
        \bm{g}_1 &= 2 \bm{\Omega} \times \bm{r}' + 2 \bm{f}
        \label{eqn:HR4BPg1old2} \\
        \bm{g}_2 &= 
        3 \bm{\Omega} \times \bm{r}' + \frac{3}{2} \bm{f} + \bm{h}_2
        \label{eqn:HR4BPg2} \\
        \bm{g}_3 &= 
        \bm{h}_3
        \label{eqn:HR4BPg3}
    \end{align}
    \label{eqn:HR4BPgk}
\end{subequations}

Note $\left( \bm{\Omega} \times \bm{r}' \right) \cdot \bm{r}' = 0$ and $\bm{f}$ corresponds to the acceleration in the CR3BP, so integrating $\bm{g}_1 \left( \prescript{*}{}{\bm{X}} (s + \tau), \tau_0 + \tau \right) \cdot \prescript{*}{}{\bm{r}'} (s + \tau)$ (i.e., computing the Melnikov function) will always result in a value of zero. The underlying culprit responsible for this result is the way time is scaled in the HR4BP, which is based on the value of $m$ according to the relationship defined in \cref{eqn:HR4BPNDt}. To first order, the primary effect from the perturbation capturing the effect of $P_0$ is associated with the scaling of time introduced in the formulation of the HR4BP.
From the previous discussion related to $\left( \bm{\Omega} \times \bm{r}' \right) \cdot \bm{r}' = 0$ and $\bm{f}$, only the $\bm{h}_k$ terms in \cref{eqn:HR4BPgk} have the potential to do any work. Let $\alpha = \tau_0 + \tau$ for conciseness. Expressions for $\bm{h}_2$ and $\bm{h}_3$ are provided in \cref{eqn:HR4BPhk}.
\begin{subequations}
    \begin{align}
        \bm{h}_2 \left( \bm{X}, \alpha \right)
        &= -\frac{3}{2} 
        \begin{bmatrix}
            -\cos{(2 \alpha)} & \sin{(2 \alpha)} & 0 \\
            \sin{(2 \alpha)} & \cos{(2 \alpha)} & 0 \\
            0 & 0 & 2/3
        \end{bmatrix}
        \bm{r}
        - \frac{1}{8} (1 - \mu) \mu \left[ P \right] 
        \begin{bmatrix}
            -8 \cos{\left( 2 \alpha \right)} \\
            11 \sin{\left( 2 \alpha \right)} \\
            0
        \end{bmatrix}
        \label{eqn:HR4BPhk2} \\
        \bm{h}_3 \left( \bm{X}, \alpha \right)
        &= - \frac{1}{12} (1 - \mu) \mu \left[ P \right] 
        \begin{bmatrix}
            -38 \cos{\left( 2 \alpha \right)} \\
            59 \sin{\left( 2 \alpha \right)} \\
            0
        \end{bmatrix}
        \label{eqn:HR4BPhk3} \\
        \left[ P \right] &=  \left( \frac{1 - \mu}{R_{1 - \mu,\text{C}}^3} -  \frac{\mu}{R_{\mu,\text{C}}^3} \right) \left[ I_{3 \times 3} \right] - \frac{3}{R_{1 - \mu,\text{C}}^5} \bm{R}_{1 - \mu,\text{C}} \bm{R}_{1 - \mu,\text{C}}^T + \frac{3}{R_{\mu,\text{C}}^5} \bm{R}_{\mu,\text{C}} \bm{R}_{\mu,\text{C}}^T
        \label{eqn:HR4BPhkP}
    \end{align}
    \label{eqn:HR4BPhk}
\end{subequations}

As we can ignore the contribution from $\bm{g}_1$, the perturbation $\varepsilon \bm{g}$ takes the form $m^2 \bm{g}_2 + m^3 \bm{g}_3 + \cdots$. So, by using $\varepsilon = m^2$ and ignoring higher order terms, the Melnikov function takes the form of:
\begin{equation}
    \mathcal{M} \left( s, \tau_0 \right) = \int_{0}^{a T^*}{\bm{h}_2 \left( \prescript{*}{}{\bm{X}} \left( s + \tau \right), \tau_0 + \tau \right) \cdot \prescript{*}{}{\bm{r}'} \left( s + \tau \right) \, \text{d}\tau}
    \label{eqn:MFHk}
\end{equation}

If the Melnikov function in \cref{eqn:MFHk} is identically zero for any $s \in \left[ 0, a T^* \right)$ on a specific CR3BP periodic orbit, then we will use the next order of $\bm{h}_j$ instead of $\bm{h}_2$ in \cref{eqn:MFHk}. Note that the validity of using a higher order $\bm{h}_j$ is contingent upon using $\varepsilon = m^j$ and the assumptions related to the behavior of the orbit presented in \cref{{sub:A2S1}} being valid with this new $\varepsilon$.

\subsection{Behavior of the HR4BP Melnikov Function}
\label{sub:PHMF}
To determine points on a CR3BP periodic orbit we expect to be able to continue into the HR4BP for $m > 0$, zeros of the Melnikov function need to be identified. This could be accomplished by selecting a value of $\tau_0$, discretizing the CR3BP periodic orbit into a set of different points (i.e., selecting different values of $s$), and then integrating \cref{eqn:MFHk} to determine the Melnikov function at each point. However, we have developed three propositions related to the behavior of the Melnikov function that demonstrate this procedure is unnecessary. Please refer to \cref{app:MTH4} for the relevant proofs. \cref{pro:P1} is generally valid for the generic form of the Melnikov function presented in \cref{eqn:MFmmain}. The other two, \cref{pro:P2} and \cref{pro:P3}, were derived for the specific form of \cref{eqn:MFHk} using the definition of $\bm{h}_2$ in \cref{eqn:HR4BPhk2}. However, it should be noted that these two propositions are also valid for the specific form of \cref{eqn:MFHk} if the definition of $\bm{h}_3$ in \cref{eqn:HR4BPhk3} is used instead of $\bm{h}_2$. Furthermore, we expect that there are similar propositions that will be valid for other of systems similar to the HR4BP, such as the ER3BP, the BCP, and the QBCP.

\begin{prop}
    Assume $\mathcal{M}(s,\tau_0)$ is a Melnikov function of the form shown in \cref{eqn:MFmmain}, where $s$ is an initial point on the unperturbed periodic orbit, and $\tau_0$ is the initial time corresponding to the periodic perturbation. Let $T = b T_g$ be the period of the periodic orbit in the perturbed system, where $b \in \mathbb{Z}_+$ and $T_g$ is the period of the forcing. Then, the following results hold for any $\tau_s \in \mathbb{R} / T \mathbb{Z}$:
    \begin{propenum}[leftmargin=*]
        \item $\mathcal{M}(s+\tau_s,\tau_0+\tau_s) = \mathcal{M}(s,\tau_0)$ \label{pro:P1A}
        \item $\mathcal{M}(s+\tau_s,\tau_0) = \mathcal{M}(s,\tau_0-\tau_s)$ \label{pro:P1B}
    \end{propenum}
    \label{pro:P1}
\end{prop}

\crefproppart{pro:P1}{pro:P1A} is a somewhat intuitive result. The Melnikov function does not change if an equal shift is applied to the initial point on the unperturbed orbit and the the initial value of time used in the integration.
\crefproppart{pro:P1}{pro:P1B} indicates that shifting the point along the initial unperturbed orbit produces the same Melnikov function as shifting the initial time of integration by an equal amount in the opposite direction. This result is known to apply to the original form of the Melnikov function as shown in~\cite{cit:Holmes}, so it is not surprising it applies to this Melnikov-like function as well.

\begin{prop}
    Assume $\mathcal{M}(s,\tau_0)$ is the Melnikov function defined in \cref{eqn:MFHk} using either $\bm{h}_2$ in \cref{eqn:HR4BPhk2} or $\bm{h}_3$ in \cref{eqn:HR4BPhk3}, $s$ is an initial point on the unperturbed periodic orbit, and $\tau_0$ is the initial time corresponding to the periodic perturbation. Let $T = b T_g$ be the period of the periodic orbit in the perturbed system, where $b \in \mathbb{Z}_+$ and $T_g$ is the period of the forcing. Then, the following relations hold for any $\tau_s \in \mathbb{R} / T \mathbb{Z}$:
    \begin{propenum}[leftmargin=*]
        \item $\mathcal{M} \left( s, \tau_0 + \tau_s \right) = \cos{(2 \tau_s)} \mathcal{M} \left( s, \tau_0 \right) - \sin{(2 \tau_s)} \mathcal{M} \left( s, \tau_0 - \frac{\pi}{4} \right)$ \label{pro:P2A}
        \item $\mathcal{M} \left( s + \tau_s, \tau_0 \right) = \cos{(2 \tau_s)} \mathcal{M} \left( s, \tau_0 \right) + \sin{(2 \tau_s)} \mathcal{M} \left( s + \frac{\pi}{4}, \tau_0 \right)$ \label{pro:P2B}
    \end{propenum}
    \label{pro:P2}
\end{prop}

\Cref{pro:P2} has two forms: \crefproppart{pro:P2}{pro:P2A} and \crefproppart{pro:P2}{pro:P2B}.
By \crefproppart{pro:P2}{pro:P2B}, the Melnikov function at any point can be determined from knowledge of the Melnikov function at two initial points on the unperturbed orbit, i.e., \cref{eqn:MFHk} must be integrated only twice to determine $\mathcal{M} \left( s, \tau_0 \right)$ and $\mathcal{M} \left( s + \frac{\pi}{4}, \tau_0 \right)$. By specifying a set of different $\tau_s$ values in the range $\left[0, a T^* \right)$, \crefproppart{pro:P2}{pro:P2B} can then be used to directly obtain $\mathcal{M} \left( s, \tau_0 \right)$ at any other value of $s = s + \tau_s$. \crefproppart{pro:P1}{pro:P1B} can then be used to obtain $\mathcal{M} \left( s, \tau_0 \right)$ at any other value of $\tau_0$. 

There are a couple other important implications of \cref{pro:P2}. First, for a particular value of $s$ and $\tau_0$, if $\mathcal{M} \left( s, \tau_0 \right) = 0$ then $\mathcal{M} \left( s + k \frac{\pi}{2}, \tau_0 \right) = 0$ for $k \in \mathbb{Z}$. Second, for a particular value of $s$ and $\tau_0$, if $\mathcal{M} \left( s, \tau_0 \right) = 0$ and $\mathcal{M} \left( s + \frac{\pi}{4}, \tau_0 \right) = 0$ then the Melnikov function is identically zero for any other value of $s$ and $\tau_0$ on that particular periodic orbit. If this is the case, then as discussed previously, we attempt to evaluate the Melnikov function using a higher order $\bm{h}_j$, provided the necessary assumptions are still valid for $\varepsilon = m^j$.

\begin{prop}
    Assume $\mathcal{M}(s,\tau_0)$ is the Melnikov function defined in \cref{eqn:MFHk} using either $\bm{h}_2$ in \cref{eqn:HR4BPhk2} or $\bm{h}_3$ in \cref{eqn:HR4BPhk3}. When evaluated at an initial point on a CR3BP periodic orbit with a period $T^*$ that satisfies the half-period symmetry conditions in \cref{eqn:HP21}, the Melnikov function takes the following simplified form:
    \begin{subequations}
        \begin{align}
            \mathcal{M} \left( s, \tau_0 \right) &=  2 A B
            \label{eqn:HP2MS1app} \\
            A &=
            \begin{cases}
              \sin{\left( 2 \tau_0 \right)}, & {\rm{if}}\ a = 1 \\
              0, & {\rm{otherwise}}
            \end{cases}
            \label{eqn:HP2MSAapp} \\
            B &= \int_{0}^{T^*/2}{ \left( K \left( \prescript{*}{}{\bm{X}} \left( s + \tau \right) \right) \cos{\left( 2 \tau \right)} - J \left( \prescript{*}{}{\bm{X}} \left( s + \tau \right) \right) \sin{\left( 2 \tau \right)} \right) \, {\rm{d}}\tau}
            \label{eqn:HP2MSBapp}
        \end{align}
        \label{eqn:HP2MSapp}
    \end{subequations}
    \label{pro:P3}
\end{prop}

Expressions for the terms $J$ and $K$ in \cref{eqn:HP2MSBapp} are provided in \cref{app:MTH4} in \cref{eqn:HP2hrd,eqn:HP2h3rd}.
Note all five of the conditions in \cref{eqn:HP21} must be met for all values of $\tau \in \left[0, T^*/2 \right]$ for the statement of half-period symmetry to be true. Some examples of such points are states on the CR3BP $L_1$ and $L_2$ planar Lyapunov and halo orbits that lie on the $xz$-plane.
\begin{subequations}
    \begin{align}
        \prescript{*}{}{x} \left( s + T^* - \tau \right) &= \prescript{*}{}{x} \left( s + \tau \right)
        \label{eqn:HP21a1} \\
        \prescript{*}{}{x'} \left( s + T^* - \tau \right) &= -\prescript{*}{}{x'} \left( s + \tau \right)
        \label{eqn:HP21a2} \\
        \prescript{*}{}{y} \left( s + T^* - \tau \right) &= -\prescript{*}{}{y} \left( s + \tau \right)
        \label{eqn:HP21b1} \\
        \prescript{*}{}{y'} \left( s + T^* - \tau \right) &= \prescript{*}{}{y'} \left( s + \tau \right)
        \label{eqn:HP21b2} \\
        \prescript{*}{}{z} \left( s + T^* - \tau \right) \prescript{*}{}{z'} \left( s + T^* - \tau \right) &= -\prescript{*}{}{z} \left( s + \tau \right) \prescript{*}{}{z'} \left( s + \tau \right)
        \label{eqn:HP21c}
    \end{align}
    \label{eqn:HP21}
\end{subequations}

For the Melnikov function in \cref{eqn:HP2MS1app} to be zero, $A$ and/or $B$ must be zero. Based on the form of $A$ and $B$, identifying when $A = 0$ is trivial. By inspection, if a CR3BP periodic orbit has a period $T^* = k \pi, \; k \in \mathbb{Z}^+$ (i.e., if $a = 1$), then any points on that orbit that satisfy the half-period symmetry conditions presented in \cref{eqn:HP21} are points where $A = 0$ provided $\tau_0 = k_1 \frac{\pi}{2}, \; k_1 \in \mathbb{Z}$ based on \cref{eqn:HP2MSAapp}. For example, let us consider the CR3BP $L_1$ and $L_2$ planar Lyapunov and halo periodic orbits with $T^* = \pi$. All points on these orbits that lie on the $xz$-plane are points where we expect to be able to continue a corresponding HR4BP periodic orbit family with $\tau_0 = 0$. If $T^*$ is not an integer multiple of $\pi$ (i.e., if $a > 1$), but there is at least one point on the orbit that satisfies \cref{eqn:HP21}, then the Melnikov function is identically zero when using $\bm{h}_2$. As using $\bm{h}_3$ instead of $\bm{h}_2$ will produce the same results, the Melnikov function should be recomputed with $\bm{h}_4$.

\subsection{Generating HR4BP Periodic Orbit Families}
\label{sub:GHPO}
The initial set of generated orbit families in the HR4BP were the five EM libration points and periodic orbits in the CR3BP whose periods were an integer multiple of $\pi$. Please note that for the remainder of this work, when a libration point is referenced, it will refer to an EM libration point. These orbit families were continued until either the maximum $m$ value at which the orbit of the primaries is stable was reached ($m = 0.19510486$)~\cite{cit:ScheeresHR4BP}, or the continuation algorithm failed to identify a new orbit member.
Other families were computed by identifying bifurcation points in this set of initial families. To identify potential bifurcation points, a singular value decomposition~(SVD) on a modified form of the corrections Jacobian was performed. The modified form of the corrections Jacobian~($\left[ DB \right]$) is presented in \cref{eqn:CJmod}.
\begin{equation}
    \left[ DB \right] = \left[ \left[ \Phi (T,0) \right]^{n_B} - \left[ I_{d \times d} \right] \; , \; \left[ \Psi_m (n_B T,0) \right] \right]
    \label{eqn:CJmod}
\end{equation}
where $ \left[ \Phi (T,0) \right]$ is the monodromy matrix of the periodic orbit, and $\left[ \Psi_m (\tau,0) \right] = \left. \frac{\partial \bm{X}}{\partial m} \right|_{\tau}$. For tangent bifurcations, period-doubling bifurcations, or other period-multiplying bifurcations ($n_B$T bifurcations where $n_B$ is an integer) in the original family, $n_B$ should be set to $1$, $2$, or $n_B$ respectively.
Note new families were only computed at tangent bifurcation points in this study.

As $\left[ DB \right]$ is a $6 \times 7$ matrix, one of its singular values will be 0. Let $\sigma_\alpha$ and $\sigma_\beta$ represent the smallest and second smallest non-zero singular values, respectively. Points on the HR4BP orbit families where a local minimum in either $\sigma_\alpha$ or $\sigma_\beta$ were identified were potential bifurcation points. To compute a different orbit family at these points, the vector ($\delta \bm{V}_\sigma$) corresponding to either $\sigma_\alpha$ or $\sigma_\beta$ (depending on which reached a local minimum) was used along with the initial state and $m$ value on the original family ($\bm{X}_0^A$ and $m^A$, respectively). The initial guess for an orbit on the new family is then represented as:

\begin{equation}
    V_0^B = \begin{bmatrix} \bm{X}_0^A \\ m^A \end{bmatrix} + \Delta s_0 \delta \bm{V}_\sigma
    \label{eqn:POBPIG}
\end{equation}

It is important to note that if the components in $\delta \bm{V}_\sigma$ corresponding to the initial state are aligned in some manner with the symmetries in \cref{eqn:HR4BPSym}, then it is possible that multiple new orbit families may exist near that bifurcation point. For example, if the initial state at $\tau = 0$ is $\bm{X}_0^A = \left[ x_0, 0, 0, 0, y'_0, 0 \right]^T$ and the components of $\delta \bm{V}_\sigma$ corresponding to the initial state are $\left[ 0, \delta y_0, 0, \delta x'_0, 0, 0 \right]^T$, then both the $+$ and $-$ directions should be used to potentially generate multiple new families. The new families were continued using the same pseudo-arclength continuation scheme used for the original families. The members of those periodic orbit families where $m = m_{\text{SEM}}$ are the periodic orbits in the HR4BP SEM system.
It should be noted that if a family has a clear symmetry, then additional orbits may be obtained by applying the appropriate transformation (e.g., to obtain the southern butterfly orbit equivalent to the northern butterfly orbit).

\section{Periodic Orbits} \label{s:POMR}
\begin{figure}[!ht]
    \centering
    \begin{minipage}[c]{0.5\linewidth}
        \centering
        \includegraphics[width=0.95\linewidth]{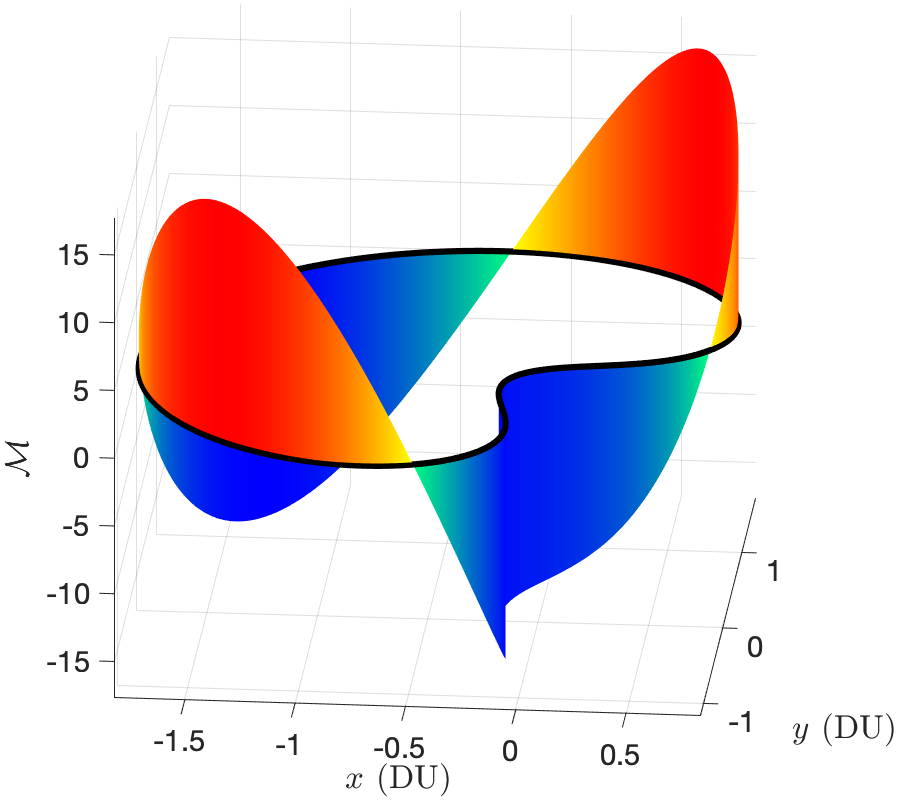}
    \end{minipage}%
    \begin{minipage}[c]{0.5\linewidth}
        \centering
        \includegraphics[width=0.95\linewidth]{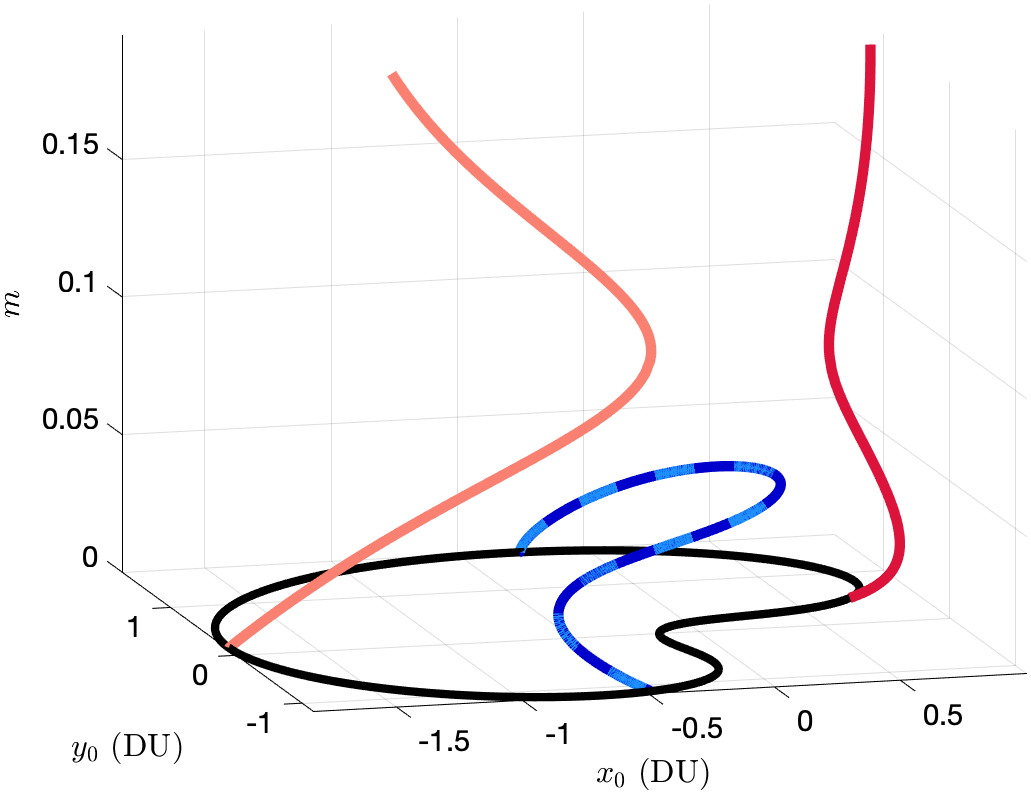}
    \end{minipage}
    \caption{Continuing the CR3BP $L_4$ Planar Orbit with $T^* = 2 \pi$ into the HR4BP. Left: Melnikov function $\mathcal{M}(s,0)$ values. Right: HR4BP periodic orbit families hodograph.}
    \label{fig:MelL42pi}
\end{figure}

The value of the Melnikov function evaluated with $\tau_0 = 0$ at different points on the CR3BP $L_4$ planar orbit member with $T^* = 2 \pi$ (so $a = 1$ and $b = 2$) is shown in \cref{fig:MelL42pi}. The four zeros of the Melnikov function are the starting points for four HR4BP orbit families which are also shown in \cref{fig:MelL42pi}.
It is important to note that these four HR4BP orbit families have been identified and computed previously by Scheeres~\cite{cit:ScheeresHR4BP} and Peterson et al.~\cite{cit:PetersonJorbaBrownEtal} using other techniques. However, the techniques we use leveraging the Melnikov function are more general. To validate our techniques, 100 points were selected on the initial CR3BP $L_4$ planar $2 \pi$ periodic orbit, and the continuation procedure was attempted starting at each of these points. Four of these points produced families in the HR4BP that could be continued up to to values of $m$ that were not negligible. These four points matched the four points where the Melnikov function is zero.

\begin{figure}[!ht]
    \centering
    \begin{minipage}[c]{0.35\linewidth}
        \centering
        \includegraphics[width=0.95\linewidth]{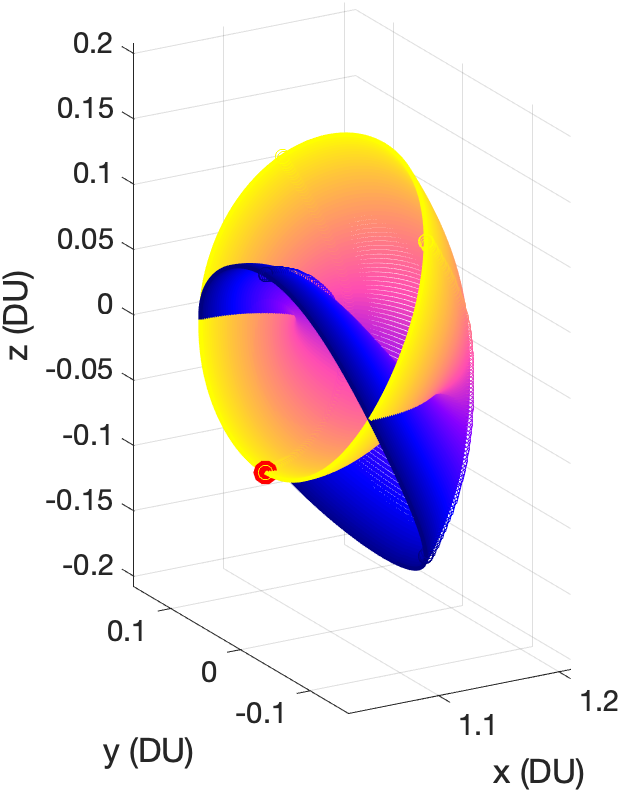}
    \end{minipage}%
    \begin{minipage}[c]{0.35\linewidth}
        \centering
        \includegraphics[width=0.95\linewidth]{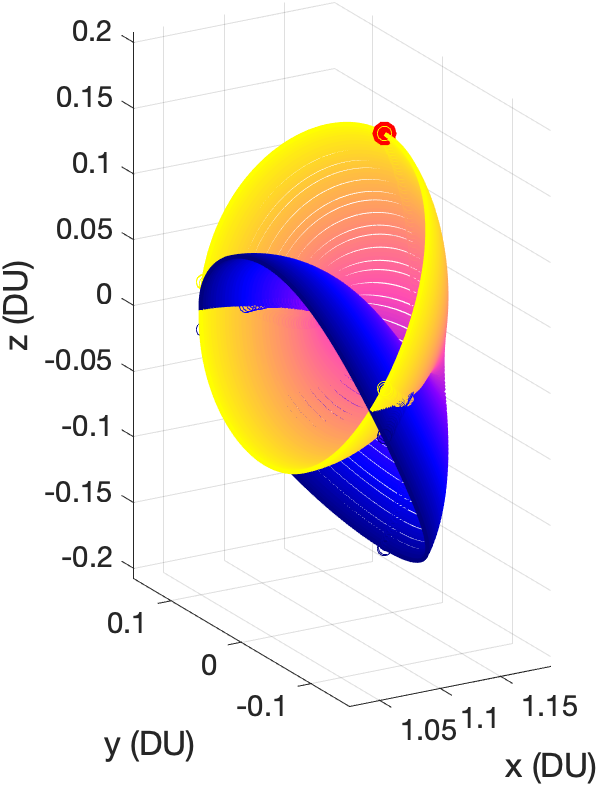}
    \end{minipage}
    \caption{HR4BP Periodic Orbit Families Corresponding to the CR3BP $L_2$ N/S Halo Orbit with $T = \pi$.}
    \label{fig:POHaloL2}
\end{figure}

Determining the number of distinct orbit families that can be generated in the HR4BP from a single orbit in the CR3BP is of interest. For the discussion related to this topic, we will distinguish between an ``orbit''/``orbit family'' and an ``object''/``object family''. An ``orbit'' refers to the collection of states on an orbit with the specific values of $\tau$ at those states.
An ``object'' refers to the collection of states on an orbit irrespective of $\tau$ at those states. Using this terminology, a single object can contain the states corresponding to multiple different orbits.
For example, two points on the CR3BP $L_2$ Northern Halo family member with $T = \pi$ can be continued (using \cref{pro:P3}) to produce two different HR4BP orbit families. At a particular value of $m$, the orbits  in these two families contain the same states, with their corresponding values of $\tau$ shifted by $\pi/2$. Therefore, only one object family in the HR4BP was identified corresponding to that particular CR3BP orbit. Note that the continuation of the orbits shown in \cref{fig:POHaloL2} starts with $m$ increasing (from yellow to magenta), before deceasing (from magenta to dark blue), until the CR3BP $L_2$ Southern Halo orbit with $T = \pi$ is obtained.

\begin{figure}[!ht]
    \centering
    \begin{minipage}[c]{0.26\linewidth}
        \centering
        \includegraphics[width=0.95\linewidth]{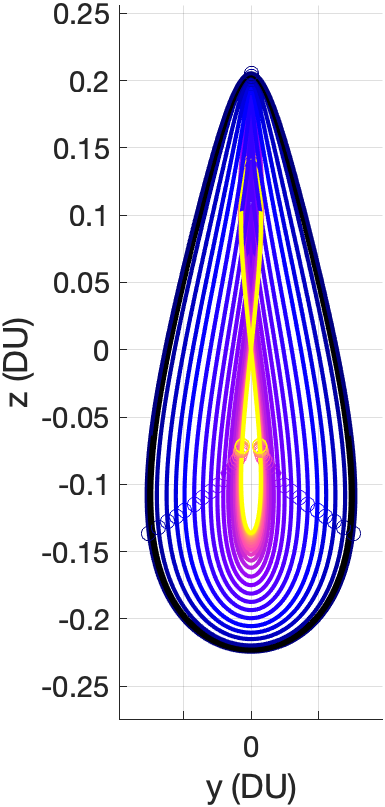}
    \end{minipage}%
    \begin{minipage}[c]{0.28\linewidth}
        \centering
        \includegraphics[width=0.95\linewidth]{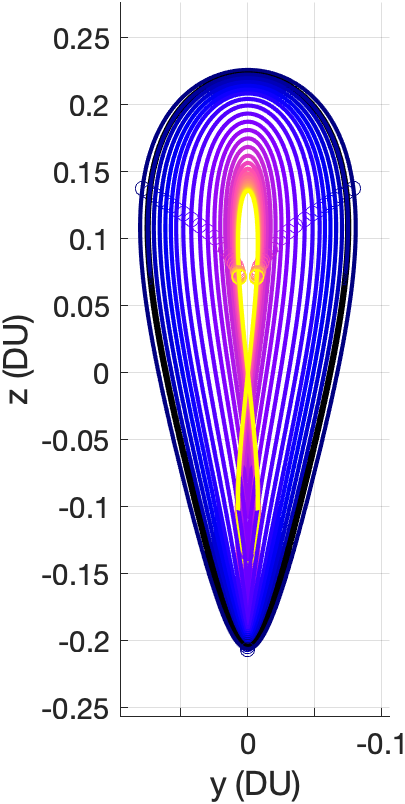}
    \end{minipage}%
    \begin{minipage}[c]{0.46\linewidth}
        \centering
        \includegraphics[width=0.95\linewidth]{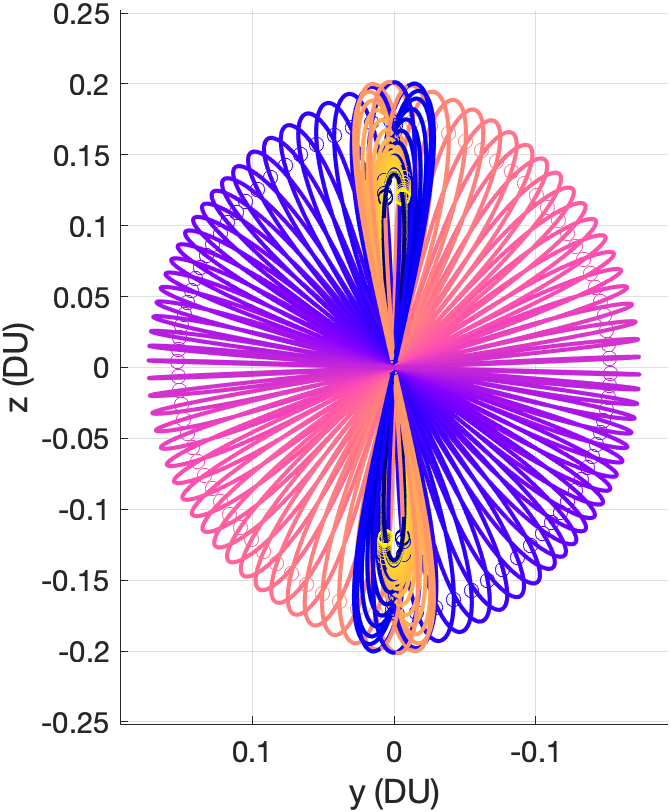}
    \end{minipage}
    \caption{HR4BP Periodic Orbit Families Corresponding to the CR3BP $L_1$ Vertical Orbit with $T^* = \pi$.}
    \label{fig:POVertL1}
\end{figure}

The $L_1$ vertical orbit with $T^* = \pi$  ($a = 1$ and $b = 1$) has four points that satisfy the symmetry conditions in \cref{eqn:HP21}. These points are separated from each other by $\tau_s = \pi/4$, so by \cref{pro:P3} and \crefproppart{pro:P2}{pro:P2B} the Melnikov function is identically zero for this orbit. We reevaluated the Melnikov function using $\bm{h}_4$ and found four zeros corresponding to the four points on the $xz$-plane. Each of these points produced a different HR4BP orbit family which belonged to one of three different object families which are shown in \cref{fig:POVertL1}.
Note that the object family on the right side of \cref{fig:POVertL1} consists of two different orbit families. The different colors are simply to distinguish between members along each family.

\subsection{Families from EM Libration Points and Associated  Periodic Orbits}

\begin{figure}[!ht]
    \centering
    \begin{minipage}[c]{0.5\linewidth}
        \centering
        \includegraphics[width=0.95\linewidth]{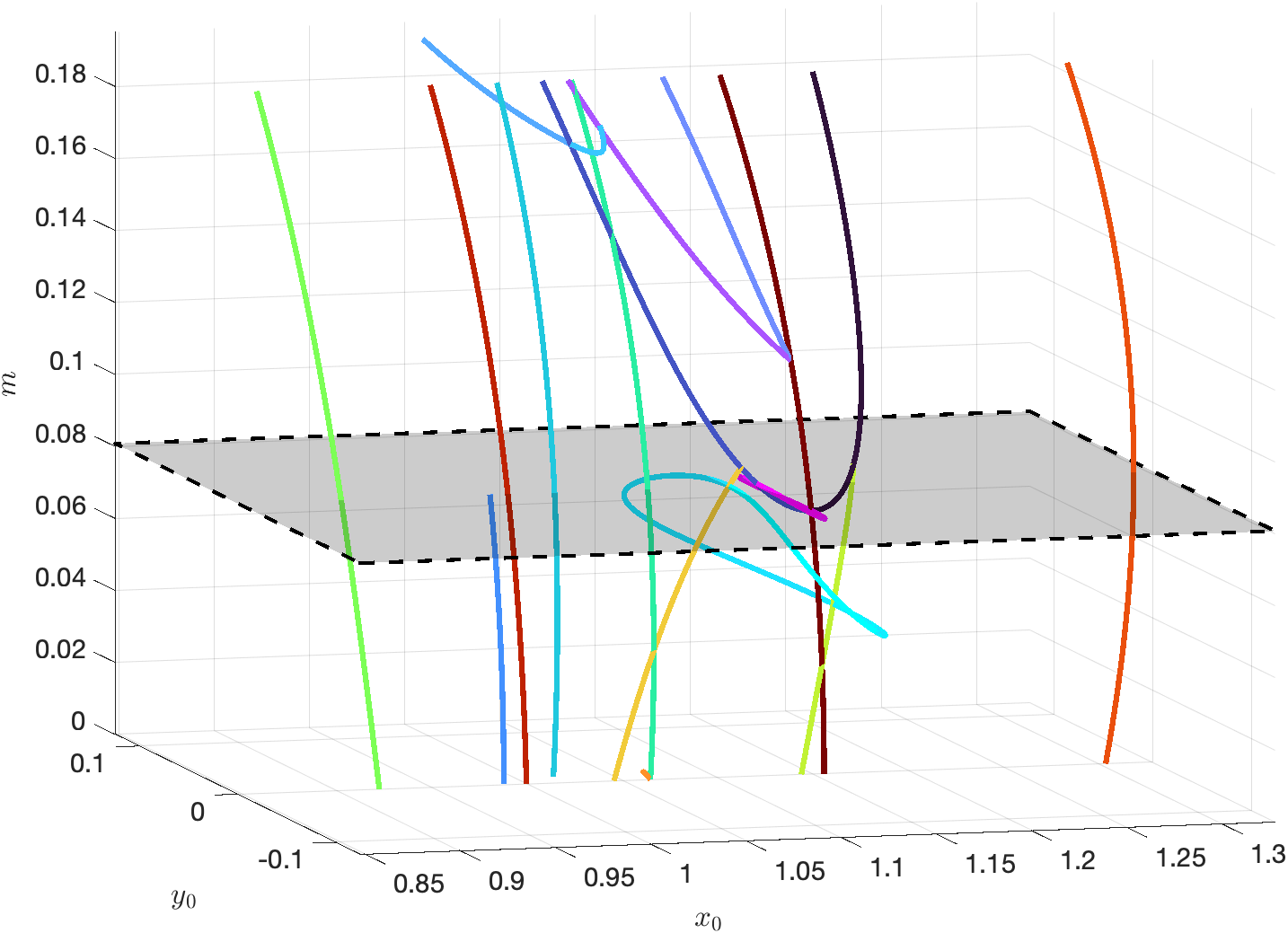}
    \end{minipage}%
    \begin{minipage}[c]{0.5\linewidth}
        \centering
        \includegraphics[width=0.95\linewidth]{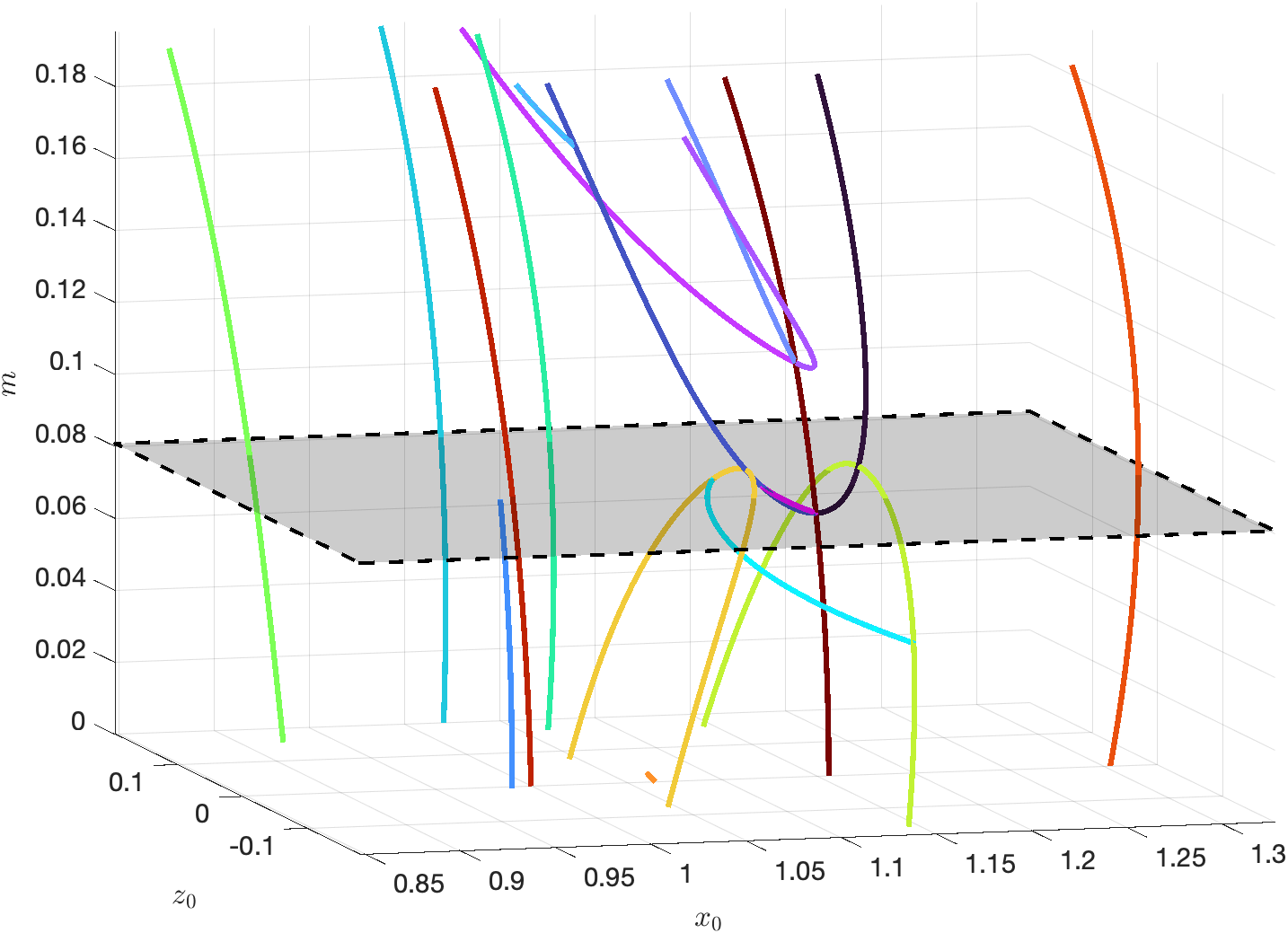}
    \end{minipage}
    \caption{Hodographs of HR4BP Periodic Orbit Families Near $L_2$ with initial states at $\tau_0 = 0$.}
    \label{fig:POBDhodL2}
\end{figure}

\begin{figure}[!ht]
    \centering
    \begin{minipage}[c]{0.5\linewidth}
        \centering
        \includegraphics[width=0.925\linewidth]{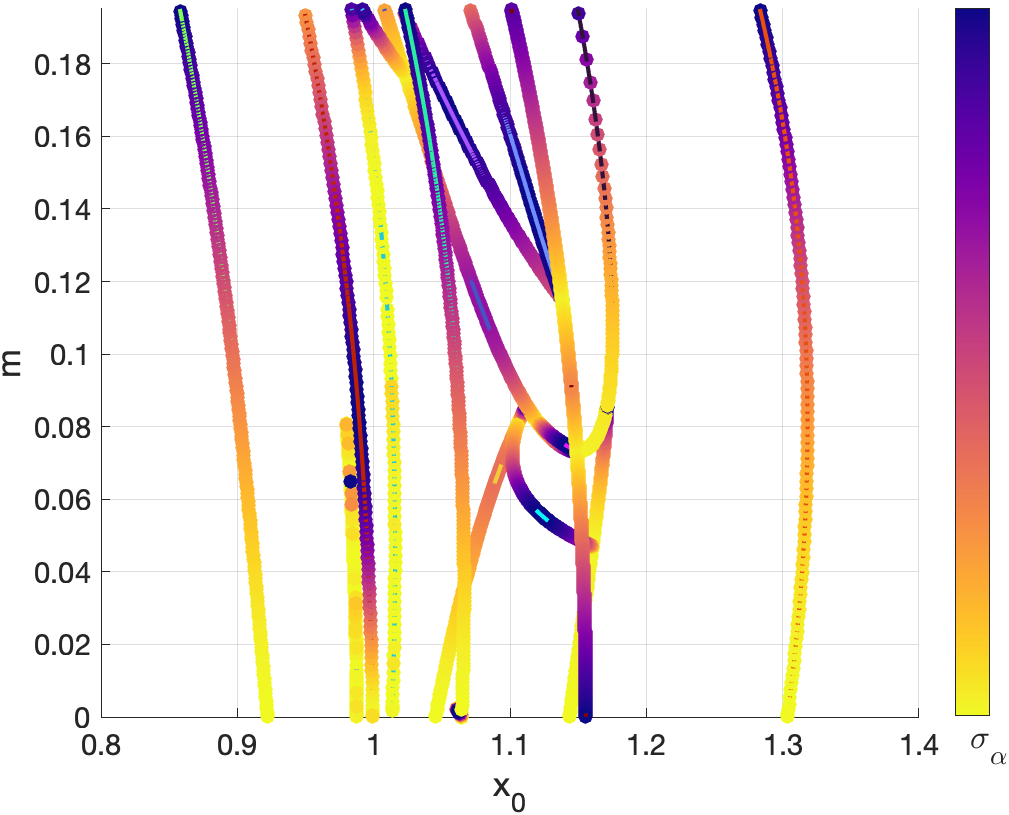}
    \end{minipage}%
    \begin{minipage}[c]{0.5\linewidth}
        \centering
        \includegraphics[width=0.925\linewidth]{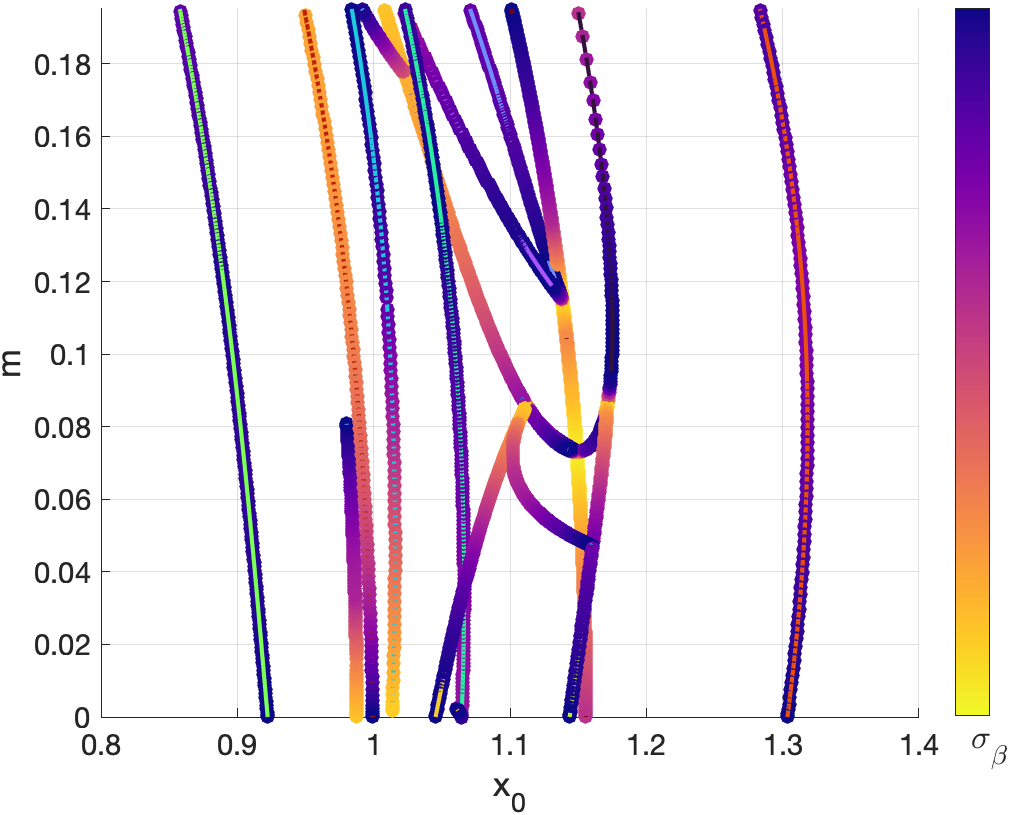}
    \end{minipage}
    \caption{$\sigma_\alpha$ (left) and $\sigma_\beta$ (right) with $n_B = 1$ for the EM $L_2$ HR4BP Periodic Orbits.}
    \label{fig:POBDsigL2}
\end{figure}

Diagrams depicting the orbit families identified around $L_2$ are presented in \cref{fig:POBDhodL2} where each color represents a different family. The gray plane represents the value of $m$ for the SEM system. The position components in these plots correspond to the initial states at $\tau_0 = 0$.
The initial set of orbits used that were related to $L_2$ were the CR3BP $L_2$ point ($T = \pi$, shown in maroon) or the A family, the planar Lyapunov family member with $T = 2 \pi$ (shown in dark red and dark orange), the vertical family member with $T = 2 \pi$ (shown in orange), the northern (and southern) halo family member with $T = \pi$ (shown in gold and light green), the northern (and southern) butterfly family member with $T = \pi$ (shown in green and aquamarine), and the 9:2 near-rectilinear halo orbit~(NRHO) with $T = 4 \pi$ (shown in teal and light blue). Families identified from bifurcations in the initial families are indicated by magenta through cyan lines.
Bifurcation diagrams showing the normalized values of $\sigma_\alpha$ and $\sigma_\beta$ along the families are provided in \cref{fig:POBDsigL2}.
Purple represents the family maximum and yellow represents zero.

Locations in \cref{fig:POBDsigL2} where local minima are reached indicate potential bifurcation points. As can be seen in these plots, many points corresponding to these local minima are where multiple families intersect.
Note that all three families associated with $L_2$ identified in~\cite{cit:OlikaraGomezMasdemont} (the A, B, and C families) were identified in this work.
Note we use a different representation of the Hill variation orbit in this work than the one used in~\cite{cit:OlikaraGomezMasdemont}. We expect this representation is the reason for the minute differences between the family structures around bifurcation points in this work and the results in~\cite{cit:OlikaraGomezMasdemont}.

\begin{figure}[!ht]
    \centering
    \begin{minipage}[c]{0.5\linewidth}
        \centering
        \includegraphics[width=0.95\linewidth]{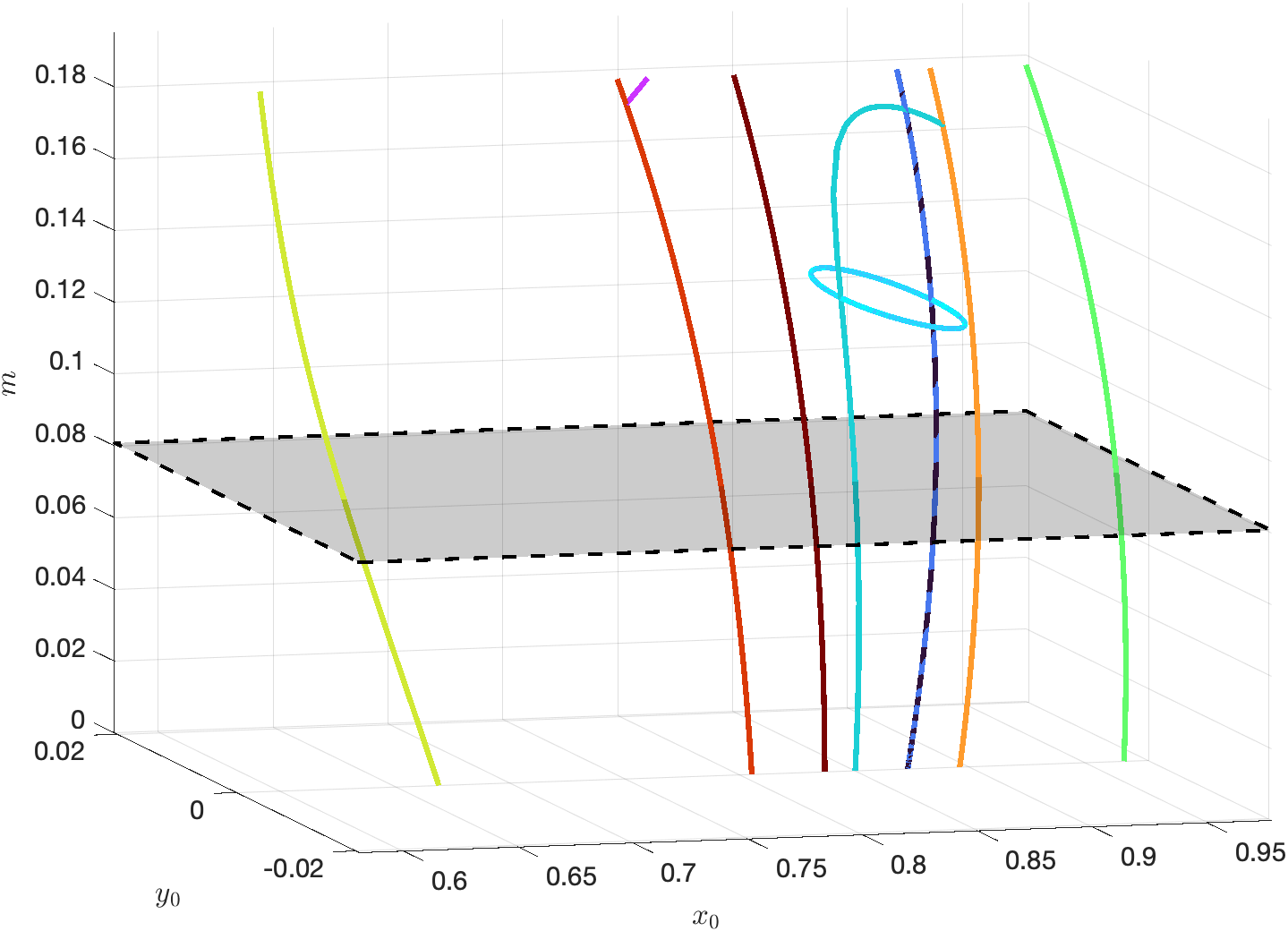}
    \end{minipage}%
    \begin{minipage}[c]{0.5\linewidth}
        \centering
        \includegraphics[width=0.95\linewidth]{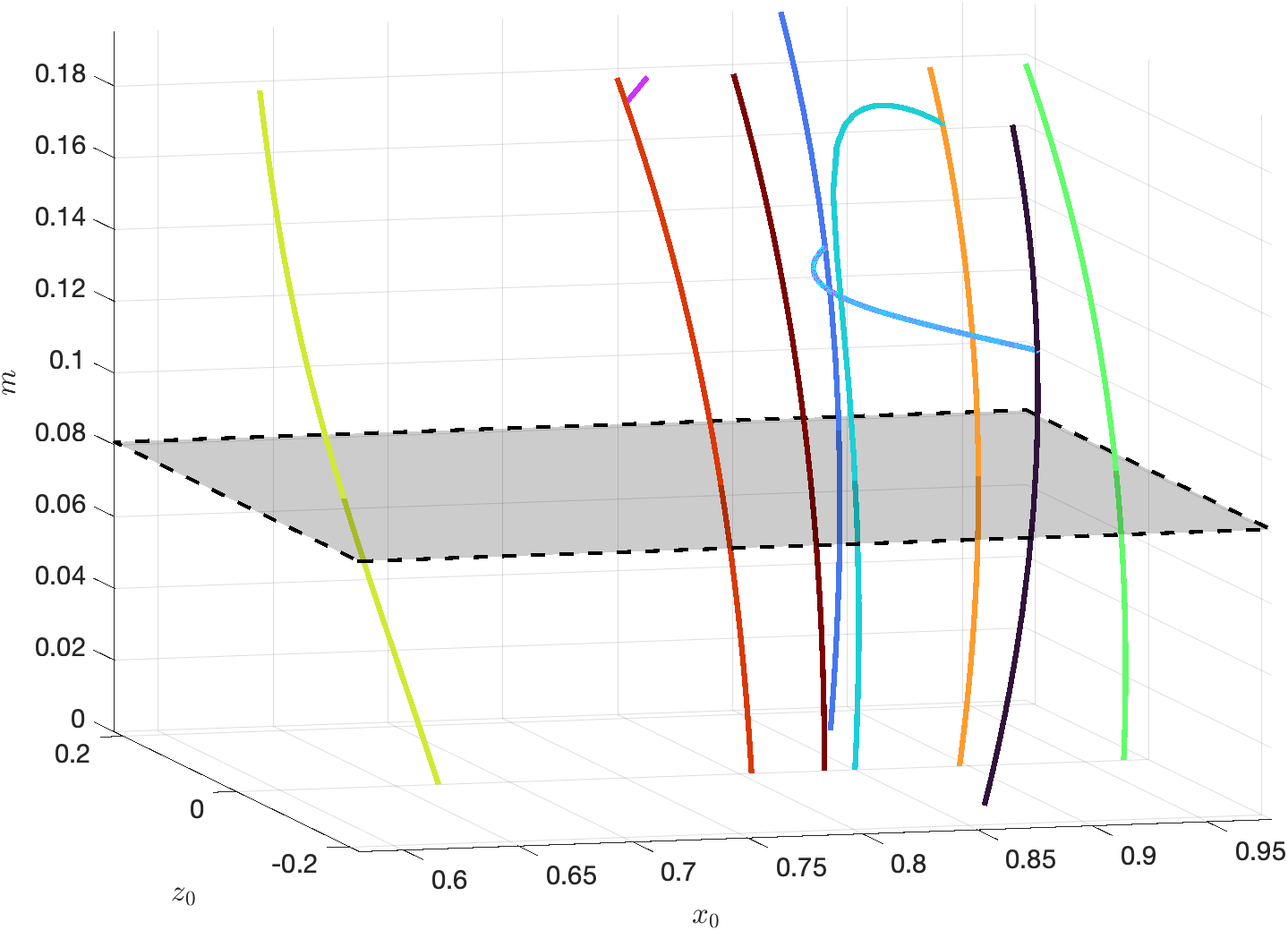}
    \end{minipage}
    \caption{Hodographs of HR4BP Periodic Orbit Families Near $L_1$ with initial states at $\tau_0 = 0$.}
    \label{fig:POBDhodL1}
\end{figure}
\begin{figure}[!ht]
    \centering
    \begin{minipage}[c]{0.5\linewidth}
        \centering
        \includegraphics[width=0.95\linewidth]{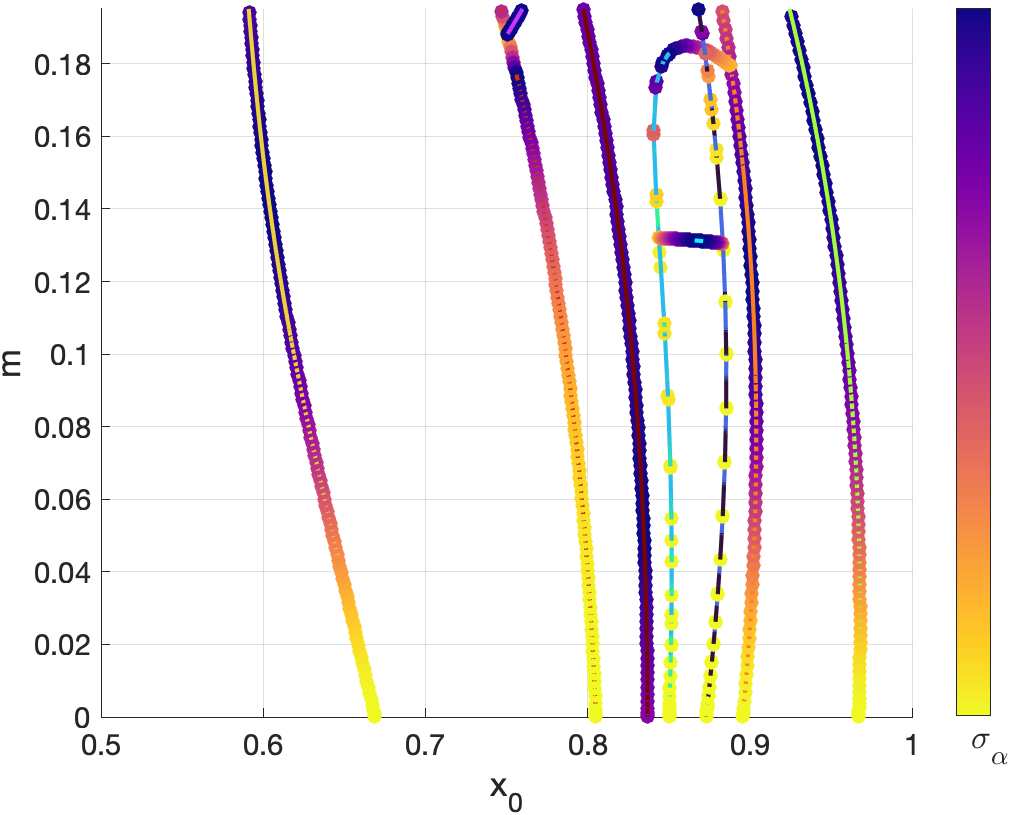}
    \end{minipage}%
    \begin{minipage}[c]{0.5\linewidth}
        \centering
        \includegraphics[width=0.95\linewidth]{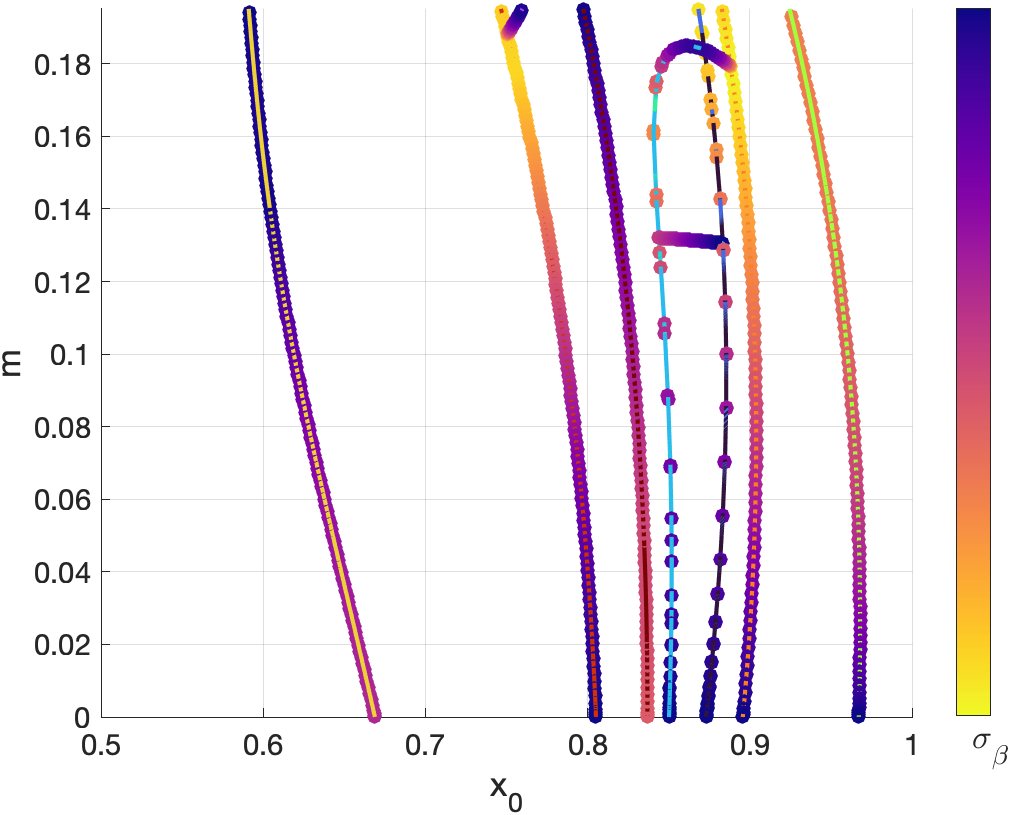}
    \end{minipage}
    \caption{$\sigma_\alpha$ (left) and $\sigma_\beta$ (right) with $n_B = 1$ for the EM $L_1$ HR4BP Periodic Orbits.}
    \label{fig:POBDsigL1}
\end{figure}

Diagrams depicting the orbit families identified around $L_1$ are presented in \cref{fig:POBDhodL1}. Bifurcation diagrams showing the normalized values of $\sigma_\alpha$ and $\sigma_\beta$ along the families are provided in \cref{fig:POBDsigL1}.
The initial set of orbits related to $L_1$ were the CR3BP $L_1$ point (with $T = \pi$, shown in maroon in \cref{fig:POBDhodL1}), the planar Lyapunov family member with $T = \pi$ (shown in red and orange) and $T = 2 \pi$ (shown in gold and green), and the vertical family member with $T = \pi$ (shown in aquamarine, light blue, and blue). Families identified from bifurcations in the initial families are indicated by the magenta through cyan lines.

\begin{figure}[!ht]
    \centering
    \begin{minipage}[c]{0.5\linewidth}
        \centering
        \includegraphics[width=0.925\linewidth]{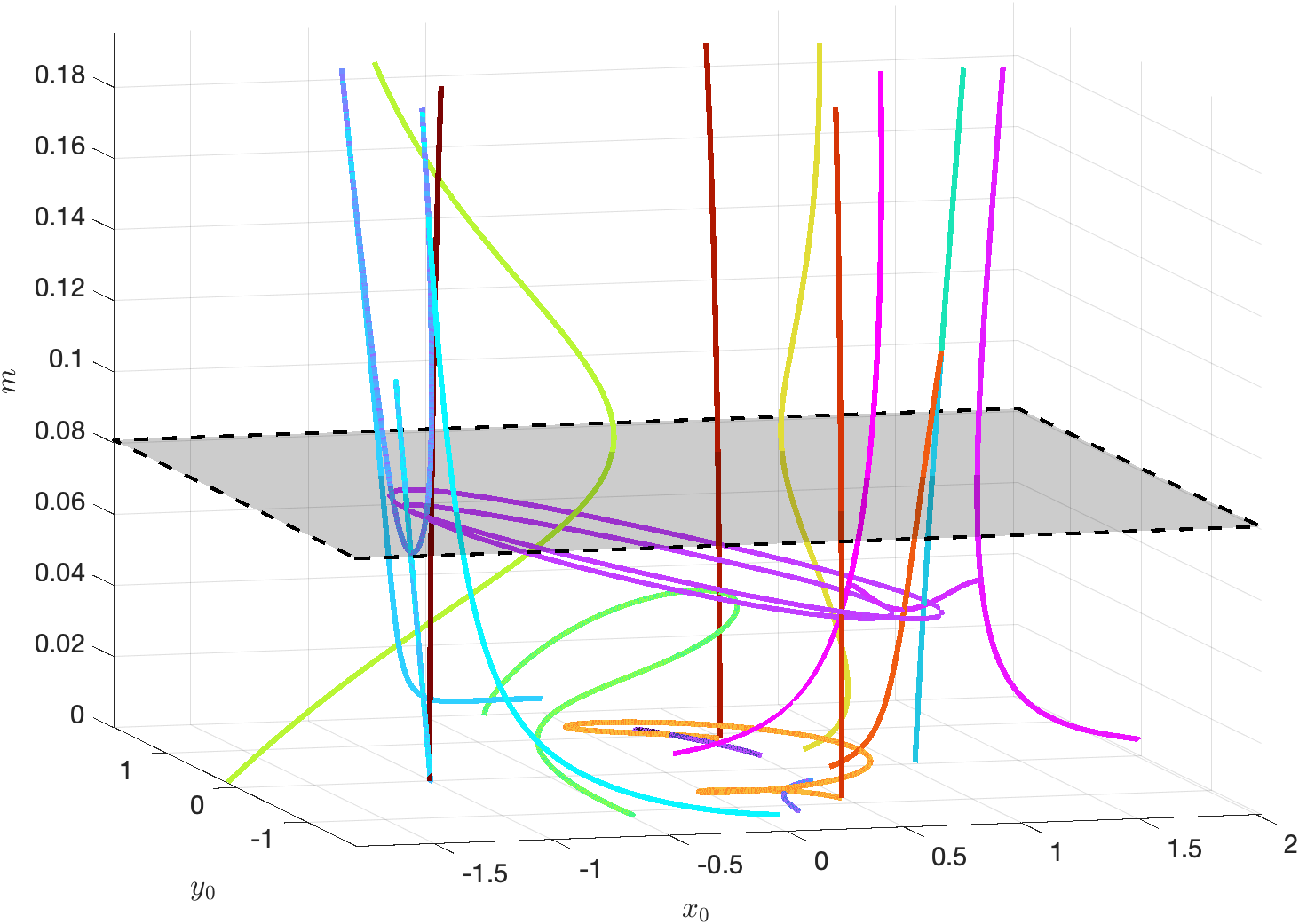}
    \end{minipage}%
    \begin{minipage}[c]{0.5\linewidth}
        \centering
        \includegraphics[width=0.925\linewidth]{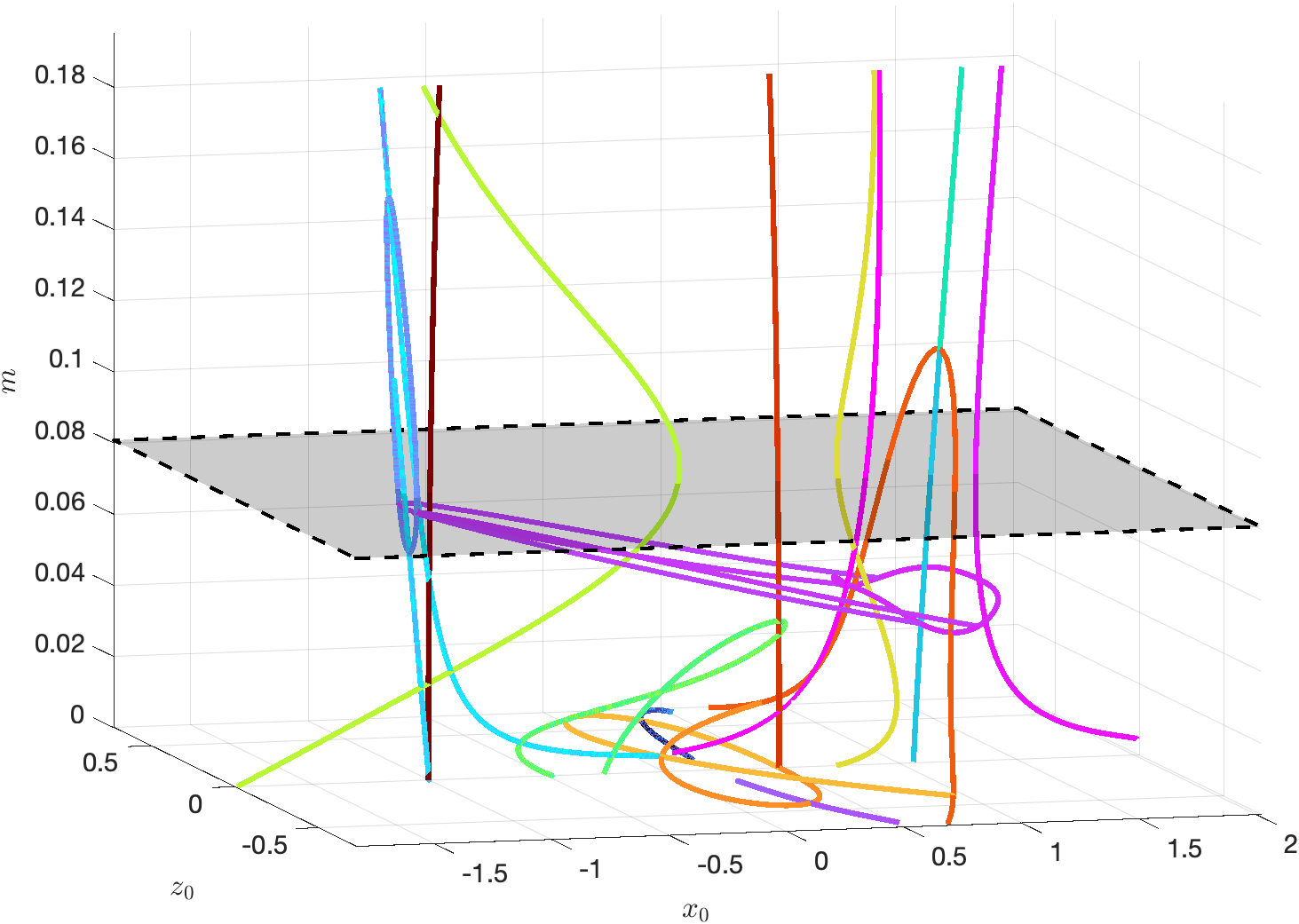}
    \end{minipage}
    \caption{Hodographs of HR4BP Periodic Orbit Families Near $L_3$, $L_4$, \& $L_5$ with states at $\tau_0 = 0$.}
    \label{fig:POBDhodL345}
    \vspace{-11pt}
\end{figure}
\begin{figure}[!ht]
    \centering
    \begin{minipage}[c]{0.5\linewidth}
        \centering
        \includegraphics[width=0.925\linewidth]{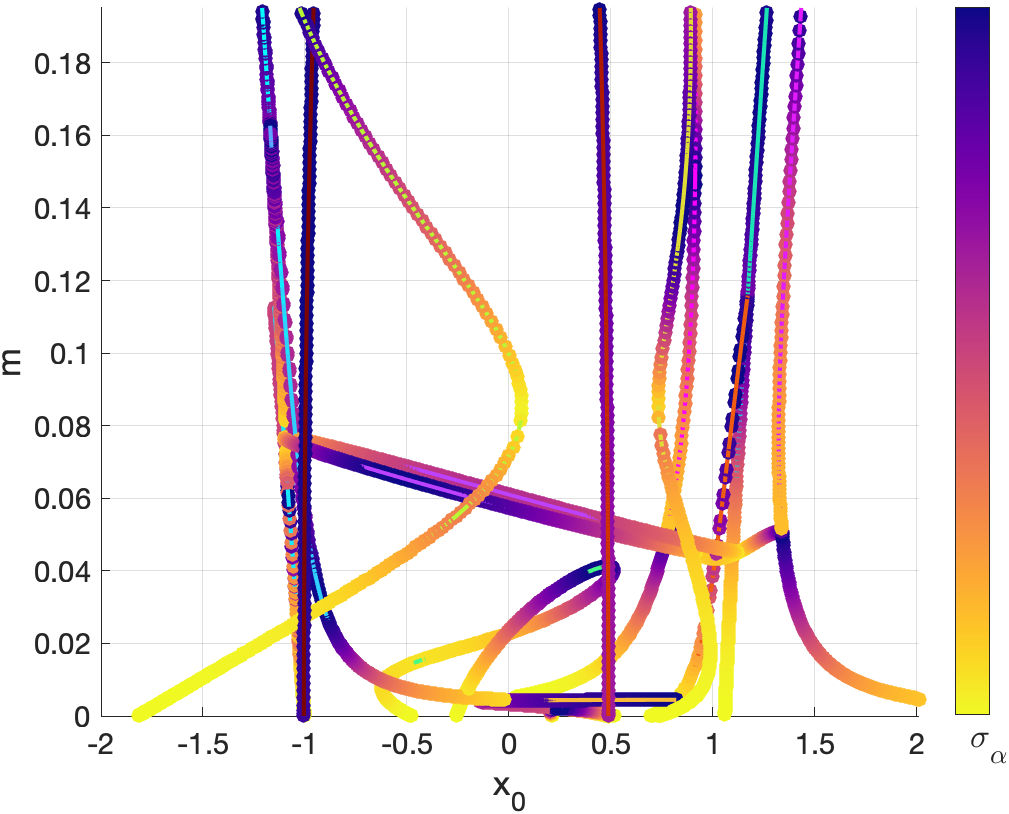}
    \end{minipage}%
    \begin{minipage}[c]{0.5\linewidth}
        \centering
        \includegraphics[width=0.925\linewidth]{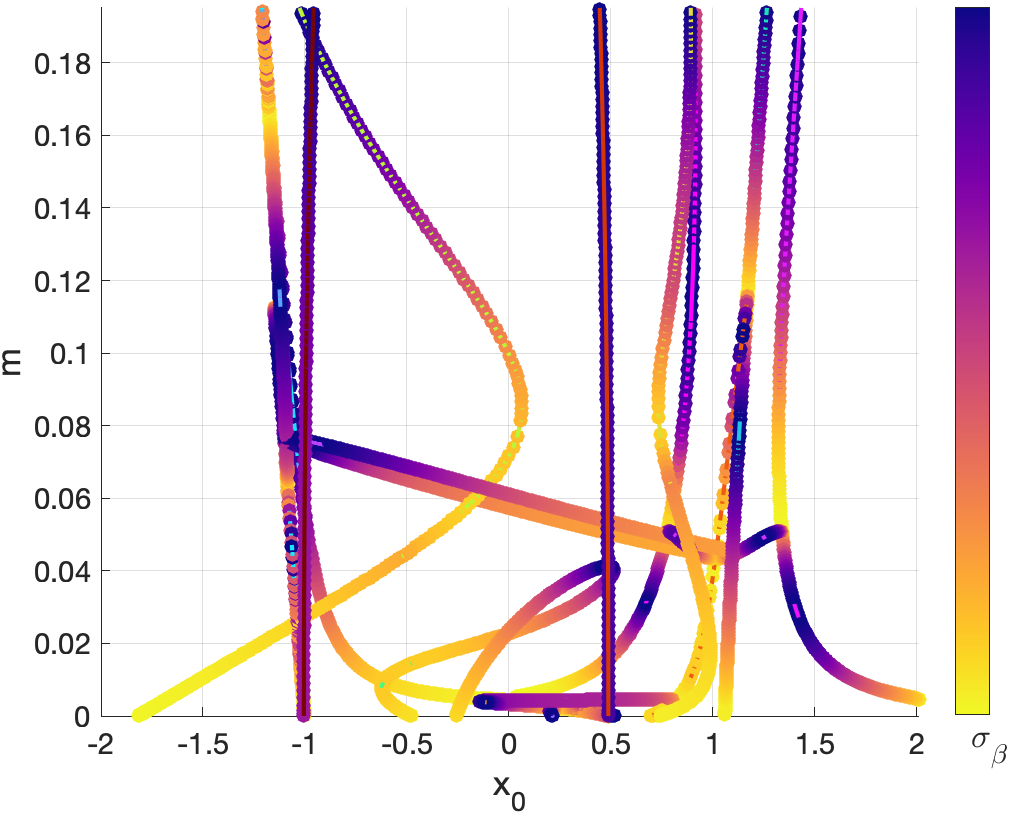}
    \end{minipage}
    \caption{$\sigma_\alpha$ (left) and $\sigma_\beta$ (right) with $n_B = 1$ for the EM $L_3$, $L_4$, \& $L_5$ HR4BP Periodic Orbits.}
    \label{fig:POBDsigL345}
    \vspace{-11pt}
\end{figure}

Diagrams depicting the orbit families identified around $L_3$, $L_4$, and $L_5$ are presented in \cref{fig:POBDhodL345}. Bifurcation diagrams showing the normalized values of $\sigma_\alpha$ and $\sigma_\beta$ along the families are provided in \cref{fig:POBDsigL345}.

The initial set of orbits related to $L_3$, $L_4$, and $L_5$ that were used were those corresponding to the CR3BP $L_3$, $L_4$, and $L_5$ points (with $T = \pi$, shown in maroon in \cref{fig:POBDhodL345}), the $L_4$ and $L_5$ planar Lyapunov family members with $T = \pi$ (shown in red and orange) and $T = 2 \pi$ (shown in gold, light green, and green), and the $L_3$, $L_4$, and $L_5$ vertical family members with $T = \pi$ (shown in aquamarine, light blue, and blue). Families identified from bifurcations in these initial families are indicated by the magenta through cyan lines.

\begin{figure}[!ht]
    \centering
    \begin{minipage}[c]{0.6\linewidth}
        \centering
        \includegraphics[width=0.95\linewidth]{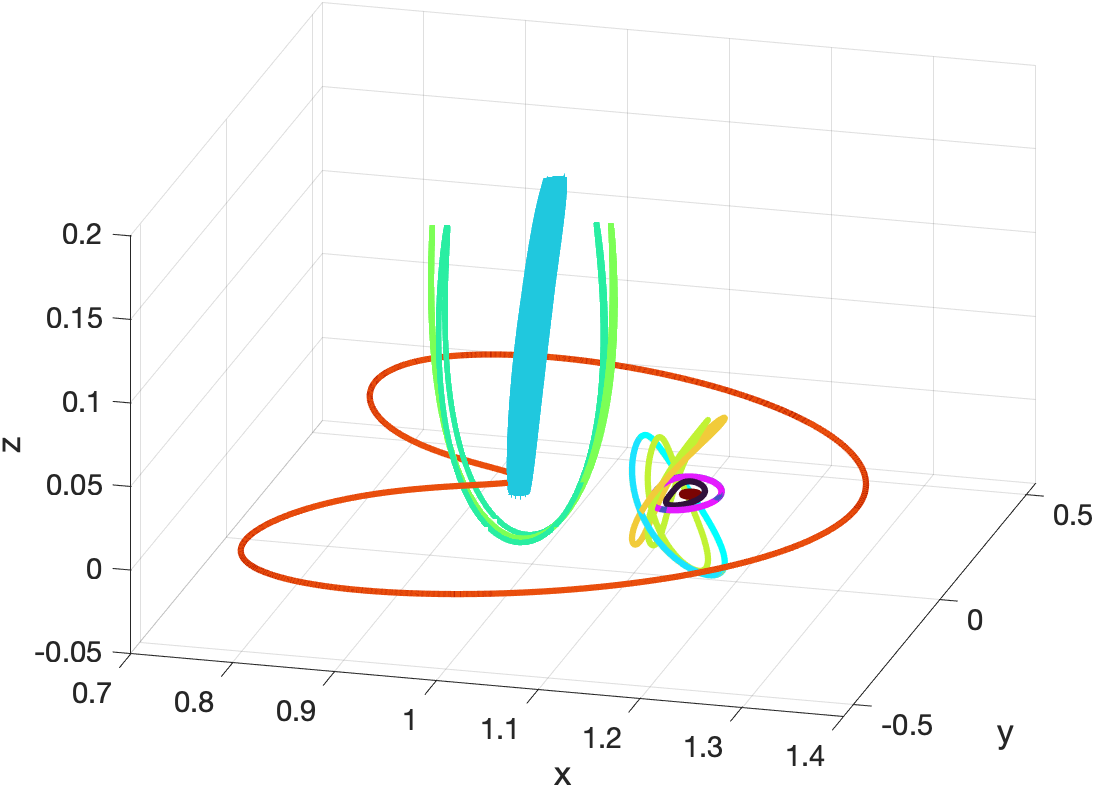}
    \end{minipage}
    \begin{minipage}[c]{0.5\linewidth}
        \centering
        \includegraphics[width=0.95\linewidth]{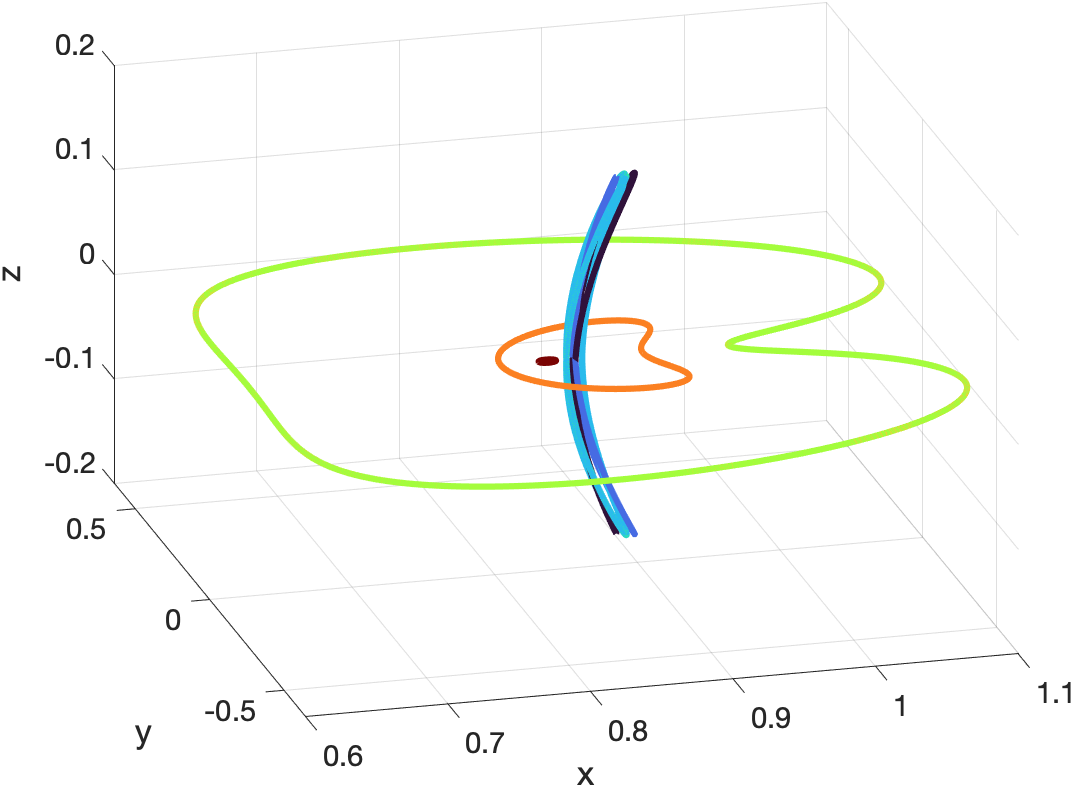}
    \end{minipage}%
    \begin{minipage}[c]{0.5\linewidth}
        \centering
        \includegraphics[width=0.95\linewidth]{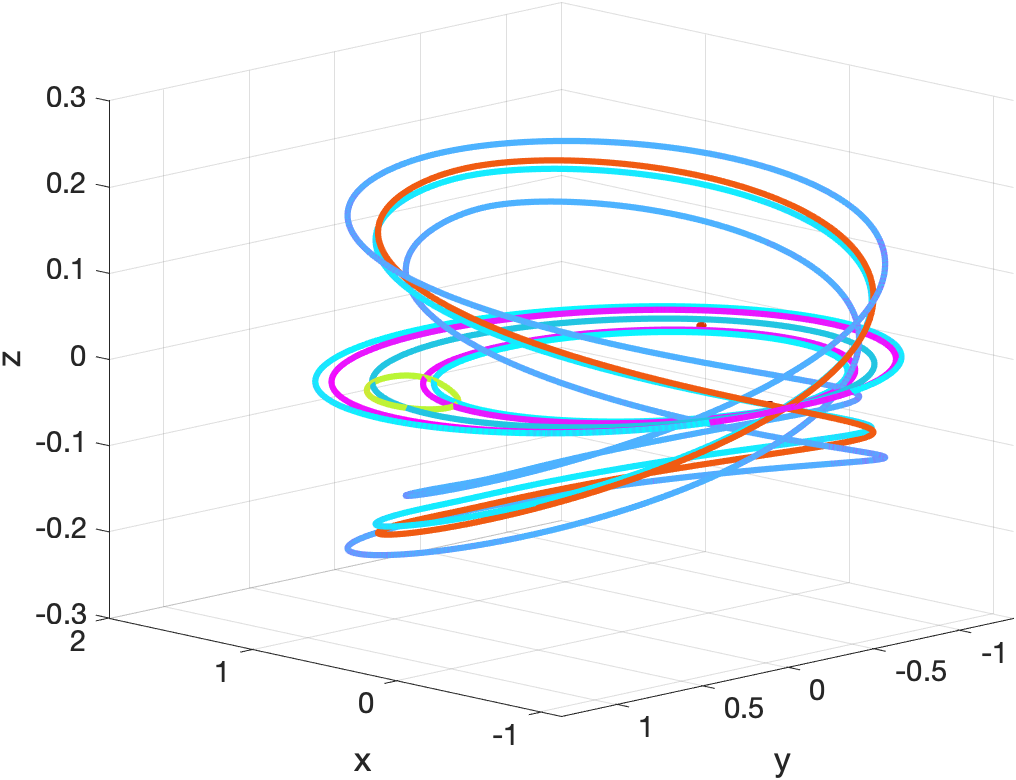}
    \end{minipage}
    \caption{HR4BP Periodic Orbits Associated with the EM Libration Points in the SEM system.}
    \label{fig:POSEM}
    \vspace{-11pt}
\end{figure}

A collection of the HR4BP periodic orbits that exist in the vicinity of EM libration points at the SEM value of $m$ are depicted in \cref{fig:POSEM}. The colors in these plots correspond to the HR4BP families in \cref{fig:POBDhodL2,fig:POBDhodL1,fig:POBDhodL345}.

\section{Conclusion} \label{s:CONC}
Periodic orbits in the CR3BP~(Circular Restricted 3-Body Problem) are useful starting points for obtaining periodic orbits in the HR4BP~(Hill Restricted 4-Body Problem). By continuing these orbits up from $m = 0$, many periodic orbit families in the HR4BP can be computed. The set of orbits in these families with $m = 0.0808$ and $\mu = 0.0122$ are periodic orbits in the HR4BP representation of the Sun-Earth-Moon system. By studying the singular values of the corrections Jacobian, leveraging symmetry when applicable, and using other continuation techniques, a variety of periodic orbits can be identified. We expect to find many other families connected to the families identified in this paper by analyzing period-multiplying bifurcations. While additional study is needed to completely map out the periodic orbit structure in the HR4BP, this work presents techniques that can be used to generate resonant periodic orbits in periodically forced systems. We validated these techniques by replicating results of previously computed orbits in the HR4BP, have extended these families to a wider parameter space, and computed additional families originating from periodic orbits around the libration points.

\clearpage

\backmatter

\bmhead{Acknowledgments}
This work was carried out with support from U.S. Air Force Office of Scientific Research grant FA9550-21-1-0332.


\appendix
\section{HR4BP Equations of Motion}
\label{app:HEOM}
The coefficients for the Hill variation orbit~(HVO) are a function of $m$ and are provided in \cref{eqn:HVO,tab:HVOdp,tab:HVOcnp}.
\begin{subequations}
    \begin{align}
        \bar{\bm{\rho}} &= \bar{\xi} \hat{\bm{\imath}}_m + \bar{\eta} \hat{\bm{\jmath}}_m + \bar{\zeta} \hat{\bm{k}}
        = \begin{bmatrix}
            \bar{\xi} \\ \bar{\eta} \\ \bar{\zeta}
        \end{bmatrix}
        = 
        \begin{bmatrix}
            \sum_{n = 1}^{N} \left( b_n + b_{-n} \right) \cos{2 n \tau} \\
            \sum_{n = 1}^{N} \left( b_{n} - b_{-n} \right) \sin{2 n \tau} \\
            0
        \end{bmatrix}
        \label{eqn:EOMS1rho} \\
        M &= \frac{m}{1 - m/3}
        \label{eqn:HVOM} \\
        a_0 &= g_0 \sum_{p = 0}^{P} d_{p} M^p
        \; \; \; \text{where} \; \; \;
        g_0 = M^{2/3}
        \label{eqn:HVOa0} \\
        b_n &= \frac{a_n}{a_0} = \sum_{p = 0}^{P} c_{n,p} M^p
        \label{eqn:HVObn}
    \end{align}
    \label{eqn:HVO}
\end{subequations}

\begin{table}[!ht]
    \caption{Coefficients $d_{p}$ for the HVO.}
    \label{tab:HVOdp}
    \centering
    \begin{tabular}{c|c}
        $p$ & $d_{p}$ \\ \hline
        0 & $1$ \\
        1 & $\frac{-8}{9}$ \\
        2 & $\frac{133}{162}$ \\
        3 & $\frac{-1264}{2187}$ \\
        4 & $\frac{3319421}{5038848}$ \\
        5 & $\frac{-13366211}{11337408}$ \\
        6 & $\frac{2028830887}{2448880128}$ \\
        7 & $\frac{-4682845907}{5509980288}$ \\
        8 & $\frac{19228022393021}{12694994583552}$ \\
        9 & $\frac{-5982128249099224247}{3119921868853739520}$
    \end{tabular}
\end{table}

\begin{table}[!ht]
    \caption{Coefficients $c_{n,p}$ for the HVO.}
    \label{tab:HVOcnp}
    \centering
    \begin{tabular}{c|cccccccc}
        & \multicolumn{8}{c}{$n$} \\
        $p$ & -4 & -3 & -2 & -1 & 1 & 2 & 3 & 4 \\ \hline
        2 & 0 & 0 & 0 & $\frac{-19}{16}$ & $\frac{3}{16}$ & 0 & 0 & 0 \\
        3 & 0 & 0 & 0 & $\frac{-7}{8}$ & $\frac{3}{8}$ & 0 & 0 & 0 \\
        4 & 0 & 0 & 0 & $\frac{11}{144}$ & $\frac{7}{48}$ & $\frac{25}{256}$ & 0 & 0 \\
        5 & 0 & 0 & $\frac{23}{640}$ & $\frac{5}{36}$ & $\frac{-1}{6}$ & $\frac{553}{1920}$ & 0 & 0 \\
        6 & 0 & $\frac{1}{192}$ & $\frac{207}{3200}$ & $\frac{-661}{82944}$ & $\frac{-34589}{110592}$ & $\frac{3743}{14400}$ & $\frac{833}{12288}$ & 0 \\
        7 & 0 & $\frac{5237}{215040}$ & $\frac{1829}{288000}$ & $\frac{374797}{276480}$ & $\frac{-22907}{46080}$ & $\frac{-28811}{864000}$ & $\frac{27337}{107520}$ & 0 \\
        8 & $\frac{23}{6144}$ & $\frac{263713}{7526400}$ & $\frac{124719}{40960000}$ & $\frac{98804551}{37324800}$ & $\frac{-23804639}{24883200}$ & $\frac{-332659139}{1105920000}$ & $\frac{5056291}{15052800}$ & $\frac{3537}{65536}$ \\
        9 & $\frac{507317}{28901376}$ & $\frac{38042489}{4741632000}$ & $\frac{48459451}{604800000}$ & $\frac{300079583}{373248000}$ & $\frac{-102469631}{124416000}$ & $\frac{-4857480211}{7257600000}$ & $\frac{472019353}{4741632000}$ & $\frac{11705987}{48168960}$
    \end{tabular}
\end{table}


These coefficients describe the Hill variation orbit family which is one particular solution to the Hill equations describing the moon's motion~\cite{cit:ScheeresHR4BP}.
The coefficients presented in \cref{tab:HVOdp,tab:HVOcnp} were derived based on the method presented by Wintner~\cite{cit:Wintner}, and are modified slightly from the form presented by Olikara and Scheeres~\cite{cit:OlikaraScheeres}.
We compute the coefficients up to a maximum order $P = 9$, so $b_n$ must be determined for integers $\left| n \right| \leq N = 4$ excluding $n = 0$. Note $c_{n,p} = 0$ for $p < 2$.

In order to evaluate the equations of motion in \cref{eqn:HR4BPEOM}, expressions for $\nabla V$, $\left[ A \right]$, and $\left[ C \right]$ are needed. Note that $\nabla V = V_{\bm{r}} = \frac{\partial V}{\partial \bm{r}}$, $\left[ \nabla \nabla V \right] = V_{\bm{r} \bm{r}} = \frac{\partial}{\partial \bm{r}} \frac{\partial V}{\partial \bm{r}}$, $V_{\bm{r} m} = \frac{\partial}{\partial m} \frac{\partial V}{\partial \bm{r}}$, and $V_{\bm{r} \mu} = \frac{\partial}{\partial \mu} \frac{\partial V}{\partial \bm{r}}$.
\begin{subequations}
    \begin{align}
        \nabla V =&
        \begin{bmatrix}
            \left( 1 + 2 m + \frac{3}{2} m^2 \right) x + \frac{3}{2} m^2 \left( x \cos{2 \tau} - y \sin{2 \tau} \right) \\
            \left( 1 + 2 m + \frac{3}{2} m^2 \right) y - \frac{3}{2} m^2 \left( y \cos{2 \tau} + x \sin{2 \tau} \right) \\
            - m^2 z 
        \end{bmatrix}
        \label{eqn:EOMD2Vr} \\
        &- \frac{m^2}{a_0^3} \left( \frac{1 - \mu}{R_{1 - \mu}^3} \bm{R}_{1 - \mu} + \frac{\mu}{R_\mu^3} \bm{R}_{\mu} \right)
        \nonumber \\
        \left[ A \right] =&
        \begin{bmatrix}
            \left[ 0_{3 \times 3} \right] & \left[ I_{3 \times 3} \right] \\
            \left[ \nabla \nabla V \right] & 2 \left(1 + m \right) \left[ S \right]
        \end{bmatrix}
        \; \text{where} \; \left[ S \right] = 
        \begin{bmatrix}
            0 & 1 & 0 \\
            -1 & 0 & 0 \\
            0 & 0 & 0
        \end{bmatrix}
        \; \text{if} \; \bm{\Omega} = \hat{\bm{k}}
        \label{eqn:EOMD2A} \\
        \left[ C \right] =&
        \begin{bmatrix}
            \left[ 0_{3 \times 1} \right] & \left[ 0_{3 \times 1} \right] \\
            -2 \bm{\Omega} \times \bm{r}' + V_{\bm{r} m}
            & V_{\bm{r} \mu}
        \end{bmatrix}
        \label{eqn:EOMD2C}
    \end{align}
    \label{eqn:EOMD2}
\end{subequations}

In order to evaluate \cref{eqn:EOMD2A,eqn:EOMD2C}, expressions for $\left[ \nabla \nabla V \right]$, $V_{\bm{r} m}$, and $V_{\bm{r} \mu}$ are needed.
\begin{subequations}
    \begin{align}
        \left[ \nabla \nabla V \right] =&
        \begin{bmatrix}
            1 + 2 m + \frac{3}{2} m^2 \left( 1 + \cos{2 \tau} \right) & -\frac{3}{2} m^2 \sin{2 \tau} & 0 \\
            -\frac{3}{2} m^2 \sin{2 \tau} & 1 + 2 m + \frac{3}{2} m^2 \left( 1 - \cos{2 \tau} \right) & 0 \\
            0 & 0 & - m^2 
        \end{bmatrix}
        \label{eqn:EOMD3Vrr} \\
        &- \frac{m^2}{a_0^3} \left( \frac{1 - \mu}{R_{1 - \mu}^3} + \frac{\mu}{R_\mu^3} \right) \left[ I_{3 \times 3} \right] + 3 \frac{m^2}{a_0^3} \frac{1 - \mu}{R_{1 - \mu}^5} \bm{R}_{1 - \mu} \bm{R}_{1 - \mu}^T + 3 \frac{m^2}{a_0^3} \frac{\mu}{R_{\mu}^5} \bm{R}_{\mu} \bm{R}_{\mu}^T
        \nonumber \\
        V_{\bm{r} m} =& 
        \begin{bmatrix}
            \left(2 + 3 m \right) x + 3 m \left( x \cos{2 \tau} - y \sin{2 \tau} \right) \\
            \left(2 + 3 m \right) y - 3 m \left( y \cos{2 \tau} + x \sin{2 \tau} \right) \\
            -2 m z
        \end{bmatrix}
        + \frac{m^2}{a_0^3} \frac{\left( 1 - \mu \right) \mu}{R_{1 - \mu}^3 R_{\mu}^3} \left( R_{1 - \mu}^3 - R_{\mu}^3 \right) \frac{\partial \bar{\bm{\rho}}}{\partial m} \label{eqn:EOMD3Vrm} \\
        &+ \frac{m}{a_0^4} \frac{1 - \mu}{R_{1 - \mu}^4} \left( 3 m a_0 \frac{\partial R_{1 - \mu}}{\partial m} - 2 a_0 R_{1 - \mu} + 3 m R_{1 - \mu} \frac{\partial a_0}{\partial m} \right) \bm{R}_{1 - \mu} \nonumber \\
        &+ \frac{m}{a_0^4} \frac{\mu}{R_{\mu}^4} \left(  3 m a_0 \frac{\partial R_{\mu}}{\partial m} - 2 a_0 R_{\mu} + 3 m R_{\mu} \frac{\partial a_0}{\partial m} \right) \bm{R}_{\mu} \nonumber \\
        V_{\bm{r} \mu} =& \frac{m^2}{a_0^3} \frac{1}{R_{1 - \mu}^4} \left( R_{1 - \mu} + 3 \left(1 - \mu \right) \frac{\partial R_{1 - \mu}}{\partial \mu} \right) \bm{R}_{1 - \mu} + \frac{m^2}{a_0^3} \frac{1}{R_{\mu}^4} \left( -R_{\mu} + 3 \mu \frac{\partial R_{\mu}}{\partial \mu} \right) \bm{R}_{\mu} \label{eqn:EOMD3Vrmu} \\
        & - \frac{m^2}{a_0^3} \left( \frac{1 - \mu}{R_{1 - \mu}^3} + \frac{\mu}{R_{\mu}^3} \right) \left( \bar{\bm{\rho}} + \hat{\bm{\imath}}_m \right) \nonumber
    \end{align}
    \label{eqn:EOMD3}
\end{subequations}

The terms needed to evaluate \cref{eqn:EOMD3} are $\frac{\partial \bar{\bm{\rho}}}{\partial m}$, $\frac{\partial a_0}{\partial m}$, $\frac{\partial R_{1 - \mu}}{\partial m}$, $\frac{\partial R_{\mu}}{\partial m}$, $\frac{\partial R_{1 - \mu}}{\partial \mu}$, and $\frac{\partial R_{\mu}}{\partial \mu}$. Expressions for all of these terms except $\frac{\partial a_0}{\partial m}$ are presented in \cref{eqn:EOMD4}.
\begin{subequations}
    \begin{align}
        \frac{\partial \bar{\bm{\rho}}}{\partial m} &= 
        \begin{bmatrix}
            \sum_{n = 1}^{N} \left( \frac{\partial b_n}{\partial m} + \frac{\partial b_{-n}}{\partial m} \right) \cos{2 n \tau} \\
            \sum_{n = 1}^{N} \left( \frac{\partial b_n}{\partial m} - \frac{\partial b_{-n}}{\partial m} \right) \sin{2 n \tau} \\
            0
        \end{bmatrix}
        \label{eqn:EOMD4drhobardm} \\
        \frac{\partial R_{1 - \mu}}{\partial m} &= \frac{\mu}{R_{1 - \mu}} \bm{R}_{1 - \mu}^T \frac{\partial \bar{\bm{\rho}}}{\partial m}
        \label{eqn:EOMD4dR1mmudm} \\
        \frac{\partial R_{\mu}}{\partial m} &= - \frac{1 - \mu}{R_{\mu}} \bm{R}_{\mu}^T \frac{\partial \bar{\bm{\rho}}}{\partial m}
        \label{eqn:EOMD4dRmudm} \\
        \frac{\partial R_{1 - \mu}}{\partial \mu} &= \frac{1}{R_{1 - \mu}} \bm{R}_{1 - \mu}^T \left( \hat{\bm{\imath}}_m + \bar{\bm{\rho}} \right)
        \label{eqn:EOMD4dR1mmudmu} \\
        \frac{\partial R_{\mu}}{\partial \mu} &= \frac{1}{R_{\mu}} \bm{R}_{\mu}^T \left( \hat{\bm{\imath}}_m + \bar{\bm{\rho}} \right)
        \label{eqn:EOMD4dRmudmu}
    \end{align}
    \label{eqn:EOMD4}
\end{subequations}

The final terms needed to evaluate the complete equations of motion are $\frac{\partial a_0}{\partial m}$ and $\frac{\partial b_n}{\partial m}$.
\begin{subequations}
    \begin{align}
        \frac{\partial a_0}{\partial m} &= \frac{\partial g_0}{\partial M} \frac{\partial M}{\partial m} \sum_{p = 0}^{P} d_{p} M^p + g_0 \frac{\partial M}{\partial m} \sum_{p = 1}^{P} p d_{p} M^{p - 1}
        \label{eqn:EOMD5a0} \\
        \frac{\partial b_n}{\partial m} &= \frac{\partial M}{\partial m} \sum_{p = 1}^{P} p c_{n,p} M^{p - 1}
        \label{eqn:EOMD5bn} \\
        \frac{\partial g_0}{\partial M} &= \frac{2}{3} M^{-1/3}
        \label{eqn:EOMD5dg0dM} \\
        \frac{\partial M}{\partial m} &= \left( 1 - \frac{m}{3} \right)^{-2}
        \label{eqn:EOMD5dMdm}
    \end{align}
    \label{eqn:EOMD5}
\end{subequations}

\section{Melnikov Theory Applications with the HR4BP}
\label{app:MTH4}
\subsection{Introduction to Melnikov Theory}
\label{sub:A2S1}
We will now discuss the derivation of the specific form of the Melnikov-like function presented by Cenedese and Haller~\cite{cit:CenedeseHaller}. Let the acceleration of an autonomous system be $\bm{r}'' = \bm{f} \left( \bm{X} \right)$, the flow of the vector field in this autonomous system be defined as $\bm{X} \left( \tau \right) = \varphi \left( \bm{X}_0, \tau \right)$, and a small time-periodic perturbation be represented by $\varepsilon \bm{g} \left( \bm{X}, \tau_0 + \tau; T_g, \varepsilon \right)$ where $\bm{g} \left( \bm{X}, \tau_0 + \tau; T_g, \varepsilon \right) = \bm{g} \left( \bm{X}, \tau_0 + \tau + T_g; T_g, \varepsilon \right)$. Let $\bm{X} \left( \tau_0 + \tau; \bm{X}_0, T_g, \varepsilon \right)$ correspond to a trajectory in the perturbed system starting at state $\bm{X}_0$ at time $\tau_0$.

Let $\bm{X}_s$ be a state on a periodic orbit in the unperturbed system with a period $T^*$ that satisfies $b T_g \approx a T^*$ where $a$ and $b$ are relatively prime integers. The state $\prescript{*}{}{\bm{X}} \left( s \right) = \varphi \left( \prescript{*}{}{\bm{X}}_s, s \right)$ for any $s \in \left[ 0, T^* \right)$ could be used as the initial state to generate a state time history of the orbit. For a particular value of $s$, the periodic orbit state time history in the unperturbed system can be represented as $\prescript{*}{}{\bm{X}} \left( s + \tau \right) = \varphi \left( \prescript{*}{}{\bm{X}}_s, s + \tau \right)$ for $\tau \in \left[ 0, T^* \right)$ where $\prescript{*}{}{\bm{X}} \left( s + \tau \right) = \prescript{*}{}{\bm{X}} \left( s + \tau + T^* \right)$ and $s \in \left[ 0, T^* \right)$.
Say we set a particular value of $s$ and $\tau_0$ and, for a small enough $\varepsilon$, there is a periodic orbit  that exists in the slightly perturbed system that corresponds to the unperturbed periodic orbit.
Let us assume that the period of this orbit is $T = b T_g = a T^* + \delta T$ where $\delta T = \mathcal{O} \left( \varepsilon \right)$. Let us also assume the initial state on this periodic orbit at $\tau_0$ is $\bm{X}_0 = \prescript{*}{}{\bm{X}} \left( s \right) + \delta \bm{X}_0$ where $\delta \bm{X}_0 = \bm{\mathcal{O}}  \left( \varepsilon \right)$. For conciseness we will represent $\bm{X} \left( \tau_0 + \tau; \bm{X}_0, T_g, \varepsilon \right)$ as $\bm{X} \left( \tau_0 + \tau; s \right)$, $\bm{g} \left( \bm{X}, \tau_0 + \tau; T_g, \varepsilon \right)$ as $\bm{g} \left( \bm{X}, \tau_0 + \tau \right)$, and $\tau_0 + \tau$ as $\alpha$. Define $\prescript{\varepsilon}{}{\delta \bm{X}} \left( \tau_0 + \tau \right) = \bm{X} \left( \tau_0 + \tau; s \right) - \prescript{*}{}{\bm{X}} \left( s + \tau \right)$ and note $\prescript{\varepsilon}{}{\delta \bm{X}} \left( 0 \right) = \delta \bm{X}_0$.
Consider the work done on the perturbed orbit by the perturbing force:
\begin{subequations}
    \begin{align}
        \delta W_{\text{P}} \left( \alpha \right) =&
        \varepsilon \bm{g} \left( \bm{X} \left( \alpha \right), \alpha \right) \cdot \bm{r}' \left( \alpha \right) \, \text{d}\tau
        \label{eqn:MTdW1}\\
        =&
        \varepsilon \left( \bm{g} \left( \prescript{*}{}{\bm{X}} \left( s + \tau \right), \alpha \right) + \left. \frac{\partial \bm{g}}{\partial \bm{X}} \right|_{\prescript{*}{}{\bm{X}} \left( s + \tau \right), \alpha} \prescript{\varepsilon}{}{\delta \bm{X}} \left( \alpha \right) + \cdots \right) \cdot \left( \prescript{*}{}{\bm{r}'} \left( s + \tau \right) + \prescript{\varepsilon}{}{\delta \bm{r}'} \left( \alpha \right) \right) \, \text{d}\tau
        \label{eqn:MTdW2}\\
        =&
        \varepsilon \bm{g} \left( \prescript{*}{}{\bm{X}} \left( s + \tau \right), \alpha \right) \cdot \prescript{*}{}{\bm{r}'} \left( s + \tau \right) \, \text{d}\tau + \varepsilon \bm{g} \left( \prescript{*}{}{\bm{X}} \left( s + \tau \right), \alpha \right) \cdot \prescript{\varepsilon}{}{\delta \bm{r}'} \left( \alpha \right) \, \text{d}\tau
        \label{eqn:MTdW3} \\
        &+ \varepsilon \left( \left. \frac{\partial \bm{g}}{\partial \bm{X}} \right|_{\prescript{*}{}{\bm{X}} \left( s + \tau \right), \alpha} \prescript{\varepsilon}{}{\delta \bm{X}} \left( \alpha \right) \right) \cdot \prescript{*}{}{\bm{r}'} \left( s + \tau \right) \, \text{d}\tau 
        \nonumber \\
        &+ \varepsilon \left( \left. \frac{\partial \bm{g}}{\partial \bm{X}} \right|_{\prescript{*}{}{\bm{X}} \left( s + \tau \right), \alpha} \prescript{\varepsilon}{}{\delta \bm{X}} \left( \alpha \right) \right) \cdot \prescript{\varepsilon}{}{\delta \bm{r}'} \left( \alpha \right) \, \text{d}\tau + \cdots
        \nonumber
    \end{align}
    \label{eqn:MTdW}
\end{subequations}

Noting the assumption $\delta \bm{X}_0 = \bm{\mathcal{O}} \left( \varepsilon \right)$, we expect the work done over the orbit related to all terms, except the first, in \cref{eqn:MTdW3} to produce results that are $\mathcal{O} \left( \varepsilon^2 \right)$ or higher.
\begin{subequations}
    \begin{align}
        W_{\text{P}} &=
        \int_{0}^{T}{\delta W_{\text{P}} \left( \alpha \right) \, \text{d}\tau}
        \label{eqn:MTW1}\\
        &=
        \varepsilon \int_{0}^{a T^* + \delta T}{\bm{g} \left( \prescript{*}{}{\bm{X}} \left( s + \tau \right), \alpha \right) \cdot \prescript{*}{}{\bm{r}'} \left( s + \tau \right) \, \text{d}\tau} + \mathcal{O} \left( \varepsilon^2 \right)
        \label{eqn:MTW2} \\
        &=
        \varepsilon \int_{0}^{a T^*}{\bm{g} \left( \prescript{*}{}{\bm{X}} \left( s + \tau \right), \tau_0 + \tau \right) \cdot \prescript{*}{}{\bm{r}'} \left( s + \tau \right) \, \text{d}\tau} + \mathcal{O} \left( \varepsilon^2 \right)
        \label{eqn:MTW3}
    \end{align}
    \label{eqn:MTW}
\end{subequations}

The previous equation is the Taylor series expansion of the so-called energy function. In this work, we refer to the leading order term of the Taylor series expansion of the energy function as the Melnikov function~\cite{cit:CenedeseHaller}.
Simplified expressions for the Melnikov function are provided in the following equation.
\begin{subequations}
    \begin{align}
        \mathcal{M} \left( s, \tau_0 \right) &= \int_{0}^{a T^*}{\bm{g} \left( \prescript{*}{}{\bm{X}} \left( s + \tau \right), \tau_0 + \tau \right) \cdot \prescript{*}{}{\bm{r}'} \left( s + \tau \right) \, \text{d}\tau}
        \label{eqn:MF1} \\
        &= \sum_{k = 0}^{a - 1}{\left[ \int_{0}^{T^*}{\bm{g} \left( \prescript{*}{}{\bm{X}} \left( s + \tau \right), \tau_0 + k T^* + \tau \right) \cdot \prescript{*}{}{\bm{r}'} \left( s + \tau \right) \, \text{d}\tau} \right]}
        \label{eqn:MF2}
    \end{align}
    \label{eqn:MF}
\end{subequations}

For an orbit to be a periodic solution, there must be no work done by the non-conservative forces on the orbit over one period. The zeros of the Melnikov function represent where this is the case when considering the leading order terms.
We expect (resonant) periodic orbits from the unperturbed system will continue into the perturbed system at points where $\mathcal{M} \left( s, \tau_0 \right) = 0$, where $\tau_0$ is the initial time of integration in the perturbed system. In the case where $\mathcal{M} \left( s, \tau_0 \right) \equiv 0$ identically, a higher-order analysis is required.
Before we present proofs of the three propositions, it is useful to note that $\prescript{*}{}{\bm{X}} \left( s + \tau + k T^* \right) = \prescript{*}{}{\bm{X}} \left( s + \tau \right)$ for any integer $k$ (i.e., $k \in \mathbb{Z}$) and $\bm{g} \left( \prescript{*}{}{\bm{X}} \left( s + \tau + a T^* \right), \tau_0 + \tau + a T^* \right) = \bm{g} \left( \prescript{*}{}{\bm{X}} \left( s + \tau \right), \tau_0 + \tau \right)$.

\begin{proof}[Proof of \cref{pro:P1}]
These results can be derived explicitly from computations. Beginning with the left-hand side of \crefproppart{pro:P1}{pro:P1A}, we separate the integral into two terms using linearity:
\begin{subequations}
    \begin{align}
        \mathcal{M} \left( s + \tau_s, \tau_0 + \tau_s \right) =& \sum_{k = 0}^{a - 1}{\left[ \int_{0}^{T^*}{\bm{g} \left( \prescript{*}{}{\bm{X}} \left( s + \tau_s + \tau \right), \tau_0 + \tau_s + k T^* + \tau \right) \cdot \prescript{*}{}{\bm{r}'} \left( s + \tau_s + \tau \right) \, \text{d}\tau} \right]}
        \label{eqn:MP11} \\
        =& \sum_{k = 0}^{a - 1}{A_k} + \sum_{k = 0}^{a - 1}{B_k}
        \label{eqn:MP13}
    \end{align}
    \label{eqn:MP1}
\end{subequations}
\noindent We then analyze the two terms independently. Note that the product of $\prescript{*}{}{\bm{r}'}$ and $\bm{g}$ (the integrand of \cref{eqn:MF1}) is periodic with a period of $a T^*$. Hence, by periodicity, we shift the integration bounds:
\begin{subequations}
    \begin{align}
        \sum_{k = 0}^{a - 1}{A_k} &= \sum_{k = 0}^{a - 1}{\left[ \int_{0}^{T^* - \tau_s}{\bm{g} \left( \prescript{*}{}{\bm{X}} \left( s + \tau_s + \tau \right), \tau_0 + \tau_s + k T^* + \tau \right) \cdot \prescript{*}{}{\bm{r}'} \left( s + \tau_s + \tau \right) \, \text{d}\tau} \right]}
        \label{eqn:MP1A1} \\
        &= \sum_{k = 0}^{a - 1}{\left[ \int_{\tau_s}^{T^*}{\bm{g} \left( \prescript{*}{}{\bm{X}} \left( s + \tau \right), \tau_0 + k T^* + \tau \right) \cdot \prescript{*}{}{\bm{r}'} \left( s + \tau \right) \, \text{d}\tau} \right]}
        \label{eqn:MP1A2}
    \end{align}
    \label{eqn:MP1A}
\end{subequations}
\noindent Applying a similar trick to the summation bounds of the second term, we obtain:
\begin{subequations}
    \begin{align}
        \sum_{k = 0}^{a - 1}{B_k} =& \sum_{k = 0}^{a - 1}{\left[ \int_{T^* - \tau_s}^{T^*}{\bm{g} \left( \prescript{*}{}{\bm{X}} \left( s + \tau_s + \tau \right), \tau_0 + \tau_s + k T^* + \tau \right) \cdot \prescript{*}{}{\bm{r}'} \left( s + \tau_s + \tau \right) \, \text{d}\tau} \right]}
        \label{eqn:MP1B1} \\
        =& \sum_{k = 0}^{a - 1}{\left[ \int_{0}^{\tau_s}{\bm{g} \left( \prescript{*}{}{\bm{X}} \left( s + \tau \right), \tau_0 + (k + 1) T^* + \tau \right) \cdot \prescript{*}{}{\bm{r}'} \left( s + \tau \right) \, \text{d}\tau} \right]}
        \label{eqn:MP1B2} \\
        =& \sum_{k = 1}^{a}{\left[ \int_{0}^{\tau_s}{\bm{g} \left( \prescript{*}{}{\bm{X}} \left( s + \tau \right), \tau_0 + k T^* + \tau \right) \cdot \prescript{*}{}{\bm{r}'} \left( s + \tau \right) \, \text{d}\tau} \right]}
        \label{eqn:MP1B3} \\
        =& \sum_{k = 1}^{a - 1}{\left[ \int_{0}^{\tau_s}{\bm{g} \left( \prescript{*}{}{\bm{X}} \left( s + \tau \right), \tau_0 + k T^* + \tau \right) \cdot \prescript{*}{}{\bm{r}'} \left( s + \tau \right) \, \text{d}\tau} \right]}
        \label{eqn:MP1B4} \\
        &+ \left[ \int_{0}^{\tau_s}{\bm{g} \left( \prescript{*}{}{\bm{X}} \left( s + \tau \right), \tau_0 + a T^* + \tau \right) \cdot \prescript{*}{}{\bm{r}'} \left( s + \tau \right) \, \text{d}\tau} \right]
        \nonumber \\
        =& \sum_{k = 1}^{a - 1}{\left[ \int_{0}^{\tau_s}{\bm{g} \left( \prescript{*}{}{\bm{X}} \left( s + \tau \right), \tau_0 + k T^* + \tau \right) \cdot \prescript{*}{}{\bm{r}'} \left( s + \tau \right) \, \text{d}\tau} \right]}
        \label{eqn:MP1B5} \\
        &+ \left[ \int_{0}^{\tau_s}{\bm{g} \left( \prescript{*}{}{\bm{X}} \left( s + \tau \right), \tau_0 + \tau \right) \cdot \prescript{*}{}{\bm{r}'} \left( s + \tau \right) \, \text{d}\tau} \right]
        \nonumber \\
        =& \sum_{k = 0}^{a - 1}{\left[ \int_{0}^{\tau_s}{\bm{g} \left( \prescript{*}{}{\bm{X}} \left( s + \tau \right), \tau_0 + k T^* + \tau \right) \cdot \prescript{*}{}{\bm{r}'} \left( s + \tau \right) \, \text{d}\tau} \right]}
        \label{eqn:MP1B6}
    \end{align}
    \label{eqn:MP1B}
\end{subequations}
\noindent Substituting these simplified expressions back into the summation for $\mathcal{M}(s+\tau_s,\tau_0+\tau_s)$ shows:
\begin{subequations}
    \begin{align}
        \mathcal{M} \left( s + \tau_s, \tau_0 + \tau_s \right) =& \sum_{k = 0}^{a - 1}{\left[ \int_{\tau_s}^{T^*}{\bm{g} \left( \prescript{*}{}{\bm{X}} \left( s + \tau \right), \tau_0 + k T^* + \tau \right) \cdot \prescript{*}{}{\bm{r}'} \left( s + \tau \right) \, \text{d}\tau} \right]}
        \label{eqn:MP1F2} \\
        &+ \sum_{k = 0}^{a - 1}{\left[ \int_{0}^{\tau_s}{\bm{g} \left( \prescript{*}{}{\bm{X}} \left( s + \tau \right), \tau_0 + k T^* + \tau \right) \cdot \prescript{*}{}{\bm{r}'} \left( s + \tau \right) \, \text{d}\tau} \right]}
        \nonumber \\
        =& \sum_{k = 0}^{a - 1}{\left[ \int_{0}^{T^*}{\bm{g} \left( \prescript{*}{}{\bm{X}} \left( s + \tau \right), \tau_0 + k T^* + \tau \right) \cdot \prescript{*}{}{\bm{r}'} \left( s + \tau \right) \, \text{d}\tau} \right]}
        \label{eqn:MP1F3} \\
        =& \mathcal{M} \left( s, \tau_0 \right)
        \label{eqn:MP1F4}
    \end{align}
    \label{eqn:MP1F}
\end{subequations}
Hence, \crefproppart{pro:P1}{pro:P1A} is proved. To prove \crefproppart{pro:P1}{pro:P1B}, we again use the periodicity of $\bm{g}$, shifting the integration bounds and separating into two terms by linearity.
\begin{subequations}
    \begin{align}
        \mathcal{M} \left( s + \tau_s, \tau_0 \right) =& \int_{0}^{a T^*}{\bm{g} \left( \prescript{*}{}{\bm{X}} \left( s + \tau_s + \tau \right), \tau_0 + \tau \right) \cdot \prescript{*}{}{\bm{r}'} \left( s + \tau_s + \tau \right) \, \text{d}\tau}
        \label{eqn:MP21} \\
        =& \int_{\tau_s}^{a T^* + \tau_s}{\bm{g} \left( \prescript{*}{}{\bm{X}} \left( s + \tau \right), \tau_0 - \tau_s + \tau \right) \cdot \prescript{*}{}{\bm{r}'} \left( s + \tau \right) \, \text{d}\tau}
        \label{eqn:MP22} \\
        =& \int_{\tau_s}^{a T^*}{\bm{g} \left( \prescript{*}{}{\bm{X}} \left( s + \tau \right), \tau_0 - \tau_s + \tau \right) \cdot \prescript{*}{}{\bm{r}'} \left( s + \tau \right) \, \text{d}\tau}
        \label{eqn:MP24} \\
        &+ \int_{0}^{\tau_s}{\bm{g} \left( \prescript{*}{}{\bm{X}} \left( s + \tau \right), \tau_0 - \tau_s + \tau \right) \cdot \prescript{*}{}{\bm{r}'} \left( s + \tau \right) \, \text{d}\tau}
        \nonumber \\
        =& \int_{0}^{a T^*}{\bm{g} \left( \prescript{*}{}{\bm{X}} \left( s + \tau \right), \tau_0 - \tau_s + \tau \right) \cdot \prescript{*}{}{\bm{r}'} \left( s + \tau \right) \, \text{d}\tau}
        \label{eqn:MP25} \\
        =& \mathcal{M} \left( s, \tau_0 - \tau_s \right)
        \label{eqn:MP26}
    \end{align}
    \label{eqn:MP2}
\end{subequations}
Hence, \crefproppart{pro:P1}{pro:P1B} is proved.
\end{proof}




\subsection{Manipulating the HR4BP Equations of Motion}
\label{sub:A2S2}
Recall the definition of the Melnikov function in \cref{eqn:MFHk} and the terms $\bm{h}_2$ and $\left[ P \right]$ in \cref{eqn:HR4BPhk}.
Note that \cref{pro:P2} and \cref{pro:P3} were derived assuming $\bm{h}_2$ is used in \cref{eqn:MFHk}. For the following proofs, it will be useful to note that $\left[ P \right]$ does not depend on $\tau_0$ (i.e., think of $\left[ P \right]$ as $\left[ P \left( \bm{X} \left( \alpha \right) \right) \right]$ which will be replaced by $\left[ P \left( \prescript{*}{}{\bm{X}} \left( s + \tau \right) \right) \right]$ when evaluating the Melnikov function).
It will also be useful to recall the following trigonometry identities.
\begin{subequations}
    \begin{align}
        \sin{(\theta_1 + \theta_2)} &= \sin{(\theta_1)} \cos{(\theta_2)} + \cos{(\theta_1)} \sin{(\theta_2)}
        \label{eqn:trigis12} \\
        \cos{(\theta_1 + \theta_2)} &= \cos{(\theta_1)} \cos{(\theta_2)} - \sin{(\theta_1)} \sin{(\theta_2)}
        \label{eqn:trigidc12} \\
        \sin{(2 \theta)} &= \cos{\left( 2 \left( \theta - \frac{\pi}{4} \right) \right)}
        \label{eqn:trigisd} \\
        \cos{(2 \theta)} &= -\sin{\left( 2 \left( \theta - \frac{\pi}{4} \right) \right)}
        \label{eqn:trigicd}
    \end{align}
    \label{eqn:trigid}
\end{subequations}


\begin{proof}[Proof of \cref{pro:P2}]
The results follow from direct computations. We begin by simplifying the form of $\bm{h}_2$ using the identities given in \cref{eqn:trigid}:
\begin{subequations}
    \begin{align}
        \bm{h}_2 \left( \bm{X} \left( \alpha \right), \alpha + \tau_s \right)
        =& -\frac{3}{2} 
        \begin{bmatrix}
            -\cos{(2 \alpha + \tau_s)} & \sin{(2 \alpha + \tau_s)} & 0 \\
            \sin{(2 \alpha + \tau_s)} & \cos{(2 \alpha + \tau_s)} & 0 \\
            0 & 0 & \frac{2}{3}
        \end{bmatrix}
        \bm{r} \left( \alpha \right)
        \label{eqn:HP1h21} \\
        &- \frac{1}{8} (1 - \mu) \mu \left[ P \left( \bm{X} \left( \alpha \right) \right) \right] 
        \begin{bmatrix}
            -8 \cos{\left( 2 \alpha + \tau_s \right)} \\
            11 \sin{\left( 2 \alpha + \tau_s \right)} \\
            0
        \end{bmatrix}
        \nonumber \\
        =& -\cos{(2 \tau_s)} \frac{3}{2} 
        \begin{bmatrix}
            -\cos{(2 \alpha)} & \sin{(2 \alpha)} & 0 \\
            \sin{(2 \alpha)} & \cos{(2 \alpha)} & 0 \\
            0 & 0 & \frac{2}{3}
        \end{bmatrix}
        \bm{r} \left( \alpha \right) 
        \label{eqn:HP1h22} \\
        &- 
        \sin{(2 \tau_s)} \frac{3}{2} 
        \begin{bmatrix}
            \sin{(2 \alpha)} & \cos{(2 \alpha)} & 0 \\
            \cos{(2 \alpha)} & -\sin{(2 \alpha)} & 0 \\
            0 & 0 & 0
        \end{bmatrix}
        \bm{r} \left( \alpha \right)
        \nonumber \\
        &- \cos{(2 \tau_s)} \frac{1}{8} (1 - \mu) \mu \left[ P \left( \bm{X} \left( \alpha \right) \right) \right] 
        \begin{bmatrix}
            -8 \cos{\left( 2 \alpha \right)} \\
            11 \sin{\left( 2 \alpha \right)} \\
            0
        \end{bmatrix}
        \nonumber \\
        &- \sin{(2 \tau_s)} \frac{1}{8} (1 - \mu) \mu \left[ P \left( \bm{X} \left( \alpha \right) \right) \right] 
        \begin{bmatrix}
            8 \sin{\left( 2 \alpha \right)} \\
            11 \cos{\left( 2 \alpha \right)} \\
            0
        \end{bmatrix}
        \nonumber \\
        =& \cos{(2 \tau_s)} \bm{h}_2 \left( \bm{X} \left( \alpha \right), \alpha \right) + \sin{(2 \tau_s)}
        \begin{bmatrix}
            0 \\
            0 \\
            z \left( \alpha \right)
        \end{bmatrix}
        - \sin{(2 \tau_s)}
        \begin{bmatrix}
            0 \\
            0 \\
            z \left( \alpha \right)
        \end{bmatrix}
        \label{eqn:HP1h23} \\
        &- 
        \sin{(2 \tau_s)} \frac{3}{2} 
        \begin{bmatrix}
            \cos{\left( 2 \left( \alpha - \frac{\pi}{4} \right) \right)} & -\sin{\left( 2 \left( \alpha - \frac{\pi}{4} \right) \right)} & 0 \\
            -\sin{\left( 2 \left( \alpha - \frac{\pi}{4} \right) \right)} & -\cos{\left( 2 \left( \alpha - \frac{\pi}{4} \right) \right)} & 0 \\
            0 & 0 & 0
        \end{bmatrix}
        \bm{r} \left( \alpha \right)
        \nonumber \\
        &- \sin{(2 \tau_s)} \frac{1}{8} (1 - \mu) \mu \left[ P \left( \bm{X} \left( \alpha \right) \right) \right] 
        \begin{bmatrix}
            8 \cos{\left( 2 \left( \alpha - \frac{\pi}{4} \right) \right)} \\
            -11 \sin{\left( 2 \left( \alpha - \frac{\pi}{4} \right) \right)} \\
            0
        \end{bmatrix}
        \nonumber \\
        =& \cos{(2 \tau_s)} \bm{h}_2 \left( \bm{X} \left( \alpha \right), \alpha \right) - 
        \sin{(2 \tau_s)} \bm{h}_2 \left( \bm{X} \left( \alpha \right), \alpha - \frac{\pi}{4} \right)
        \label{eqn:HP1h24} \\
        &- \sin{(2 \tau_s)} z \left( \alpha \right) \hat{\bm{k}}
        \nonumber
    \end{align}
    \label{eqn:HP1h2}
\end{subequations}
Substituting the simplified form of $\bm{h}_2$ into the equation for $\mathcal{M} (s, \tau_0 + \tau_s)$, we obtain:
\begin{subequations}
    \begin{align}
        \mathcal{M} \left( s, \tau_0 + \tau_s \right) =& \int_{0}^{a T^*}{\bm{h}_2 \left( \prescript{*}{}{\bm{X}} \left( s + \tau \right), \alpha + \tau_s \right) \cdot \prescript{*}{}{\bm{r}'} \left( s + \tau \right) \, \text{d}\tau}
        \label{eqn:HP11} \\
        =& \int_{0}^{a T^*}{\cos{(2 \tau_s)} \bm{h}_2 \left( \prescript{*}{}{\bm{X}} \left( s + \tau \right), \alpha \right) \cdot \prescript{*}{}{\bm{r}'} \left( s + \tau \right) \, \text{d}\tau}
        \label{eqn:HP12} \\
        &+ \int_{0}^{a T^*}{-\sin{(2 \tau_s)} \bm{h}_2 \left( \prescript{*}{}{\bm{X}} \left( s + \tau \right), \alpha - \frac{\pi}{4} \right) \cdot \prescript{*}{}{\bm{r}'} \left( s + \tau \right) \, \text{d}\tau}
        \nonumber \\
        &+ \int_{0}^{a T^*}{-\sin{(2 \tau_s)} \prescript{*}{}{z} \left( s + \tau \right) \hat{\bm{k}}
        \cdot \prescript{*}{}{\bm{r}'} \left( s + \tau \right) \, \text{d}\tau}
        \nonumber \\
        =& \cos{(2 \tau_s)} \int_{0}^{a T^*}{\bm{h}_2 \left( \prescript{*}{}{\bm{X}} \left( s + \tau \right), \tau_0 + \tau \right) \cdot \prescript{*}{}{\bm{r}'} \left( s + \tau \right) \, \text{d}\tau}
        \label{eqn:HP13} \\
        &-\sin{(2 \tau_s)} \int_{0}^{a T^*}{\bm{h}_2 \left( \prescript{*}{}{\bm{X}} \left( s + \tau \right), \tau_0 - \frac{\pi}{4} + \tau \right) \cdot \prescript{*}{}{\bm{r}'} \left( s + \tau \right) \, \text{d}\tau}
        \nonumber \\
        &-\sin{(2 \tau_s)} \int_{0}^{a T^*}{\prescript{*}{}{z} \left( s + \tau \right)
        \cdot \prescript{*}{}{z'} \left( s + \tau \right) \, \text{d}\tau}
        \nonumber \\
        =& \cos{(2 \tau_s)} \mathcal{M} \left( s, \tau_0 \right) - \sin{(2 \tau_s)} \mathcal{M} \left( s, \tau_0 - \frac{\pi}{4} \right) - \sin{(2 \tau_s)} \left( 0 \right)
        \label{eqn:HP14} \\
        =& \cos{(2 \tau_s)} \mathcal{M} \left( s, \tau_0 \right) - \sin{(2 \tau_s)} \mathcal{M} \left( s, \tau_0 - \frac{\pi}{4} \right)
        \label{eqn:HP15}
    \end{align}
    \label{eqn:HP1}
\end{subequations}
proving \crefproppart{pro:P2}{pro:P2A}. \crefproppart{pro:P2}{pro:P2B} follows from an application of \crefproppart{pro:P1}{pro:P1B}:
\begin{subequations}
    \begin{align}
        \mathcal{M} \left( s + \tau_s, \tau_0 \right) &= \mathcal{M} \left( s, \tau_0 - \tau_s \right)
        \label{eqn:HP1B1} \\
        &= \cos{(-2 \tau_s)} \mathcal{M} \left( s, \tau_0 \right) - \sin{(-2 \tau_s)} \mathcal{M} \left( s, \tau_0 - \frac{\pi}{4} \right)
        \label{eqn:HP1B2} \\
        &= \cos{(2 \tau_s)} \mathcal{M} \left( s, \tau_0 \right) + \sin{(2 \tau_s)} \mathcal{M} \left( s, \tau_0 - \frac{\pi}{4} \right)
        \label{eqn:HP1B3} \\
        &= \cos{(2 \tau_s)} \mathcal{M} \left( s, \tau_0 \right) + \sin{(2 \tau_s)} \mathcal{M} \left( s + \frac{\pi}{4}, \tau_0 \right)
        \label{eqn:HP1B4}
    \end{align}
    \label{eqn:HP1B}
\end{subequations}
Hence, \crefproppart{pro:P2}{pro:P2B} is proved.
\end{proof}

The final proposition we will prove applies to periodic orbits in the CR3BP that exhibit a half-period symmetry condition at some $s \in \mathbb{R}$, which occurs often in the computations and analyses contained in the present article.

\begin{proof}[Proof of \cref{pro:P3}]
We will begin by obtaining a more detailed expression for the $\left[ P \right]$ matrix presented in \cref{eqn:HR4BPhkP}. For conciseness we will represent $\bm{X} \left( \alpha \right)$ as $\bm{X}$.  

\begin{subequations}
    \begin{align}
        \left[ P \right] &= 
        \begin{bmatrix}
            a & - \left( d \left( x + \mu \right) + e \right) y & - \left( d \left( x + \mu \right) + e \right) z \\
            - \left( d \left( x + \mu \right) + e \right) y & b & -d y z \\
            - \left( d \left( x + \mu \right) + e \right) z & -d y z & c
        \end{bmatrix}
        \label{eqn:HP2PP} \\
        a &= \kappa_3 - 3 \left( x + \mu \right)^2 \kappa_5 - 3 \frac{2 \left( x + \mu \right) - 1}{R_{\mu,\text{C}}^5}
        \label{eqn:HP2Pa} \\
        b &= \kappa_3 - 3 y^2 \kappa_5
        \label{eqn:HP2Pb} \\
        c &= \kappa_3 - 3 z^2 \kappa_5
        \label{eqn:HP2Pc} \\
        d &= 3 \kappa_5
        \label{eqn:HP2Pd} \\
        e &= \frac{3}{R_{\mu,\text{C}}^5}
        \label{eqn:HP2Pe} \\
        \kappa_k &= \frac{1}{R_{1 - \mu,\text{C}}^k} - \frac{1}{R_{\mu,\text{C}}^k}
        \label{eqn:HP2Pkk}
    \end{align}
    \label{eqn:HP2P}
\end{subequations}

Note the values of $a$, $b$, $c$, $d$, $e$, $\kappa_3$, and $\kappa_5$ are the same when evaluated at $\prescript{*}{}{\bm{X}} \left( s + T^* - \tau \right)$ and $\prescript{*}{}{\bm{X}} \left( s + \tau \right)$ if the point $s$ satisfies the half-period symmetry conditions presented in \cref{eqn:HP21}. We will now use the result in \cref{eqn:HP2PP} to obtain another expression for $\bm{h}_2$, which was presented in \cref{eqn:HR4BPhk2}.
\begin{equation}
    \begin{split}
        \bm{h}_2 \left( \bm{X}, \alpha \right) =& \left(
        \frac{3}{2} 
        \begin{bmatrix}
            x \\ -y \\ 0
        \end{bmatrix}
        - \left( 1 - \mu \right) \mu 
        \begin{bmatrix}
            -a \\ y \left( d \left( x + \mu \right) + e \right) \\ z \left( d \left( x + \mu \right) + e \right)
        \end{bmatrix}
        \right) \cos{(2 \alpha)}
        \\
        &+ \left(
        \frac{3}{2} 
        \begin{bmatrix}
            -y \\ -x \\ 0
        \end{bmatrix}
        + \frac{11}{8} \left( 1 - \mu \right) \mu 
        \begin{bmatrix}
            y  \left( d \left( x + \mu \right) + e \right) \\ -b \\ d y z
        \end{bmatrix}
        \right) \sin{(2 \alpha)}
        +
        \begin{bmatrix}
            0 \\ 0 \\ -z
        \end{bmatrix}
    \end{split}
    \label{eqn:HP2h2}
\end{equation}
Evaluating the dot product of $\bm{h}_2$ with the velocity $\bm{r}'$ yields:
\begin{subequations}
    \begin{align}
        \bm{h}_2 \left( \bm{X}, \alpha \right) \cdot \bm{r}' &= 
        J \left( \bm{X} \right) \cos{(2 \alpha)} + K \left( \bm{X} \right) \sin{(2 \alpha)} - z z'
        \label{eqn:HP2hrd1} \\
        J \left( \bm{X} \right) &= \frac{3}{2} \left( x x' - y y' \right) - \left( 1 - \mu \right) \mu \left( -a x' + \left( y y' + z z' \right) \left( d \left( x + \mu \right) + e \right) \right)
        \label{eqn:HP2hrdJ} \\
        K \left( \bm{X} \right) &= \frac{3}{2} \left( -y x' - x y' \right) + \frac{11}{8} \left( 1 - \mu \right) \mu \left( y x' \left( d \left( x + \mu \right) + e \right) - b y' + d y z z' \right)
        \label{eqn:HP2hrdK}
    \end{align}
    \label{eqn:HP2hrd}
\end{subequations}
Note that $J \left( \prescript{*}{}{\bm{X}} \left( s + T^* - \tau \right) \right) = - J \left( \prescript{*}{}{\bm{X}} \left( s + \tau \right) \right)$ and $K \left( \prescript{*}{}{\bm{X}} \left( s + T^* - \tau \right) \right) = K \left( \prescript{*}{}{\bm{X}} \left( s + \tau \right) \right)$ if the point corresponding to $s$ satisfies the half-period symmetry conditions presented in \cref{eqn:HP21}.
Note a $u$ substitution is performed on the second integral between \cref{eqn:HP2M2} and \cref{eqn:HP2M3} where $u = T^* - \tau$. 
The variable $\tau$ is still used instead of $u$ as the name of the variable is arbitrary.
\begin{subequations}
    \begin{align}
        \mathcal{M} \left( s, \tau_0 \right) =& \sum_{k = 0}^{a - 1}{\left[ \int_{0}^{T^*}{\bm{h}_2 \left( \prescript{*}{}{\bm{X}} \left( s + \tau \right), \tau_0 + k T^* + \tau \right) \cdot \prescript{*}{}{\bm{r}'} \left( s + \tau \right) \, \text{d}\tau} \right]}
        \label{eqn:HP2M1} \\
        =& \sum_{k = 0}^{a - 1}{\left[ \int_{0}^{T^*/2}{\bm{h}_2 \left( \prescript{*}{}{\bm{X}} \left( s + \tau \right), \tau_0 + k T^* + \tau \right) \cdot \prescript{*}{}{\bm{r}'} \left( s + \tau \right) \, \text{d}\tau} \right]}
        \label{eqn:HP2M2} \\
        &+ \sum_{k = 0}^{a - 1}{\left[ - \int_{T^*}^{T^*/2}{\bm{h}_2 \left( \prescript{*}{}{\bm{X}} \left( s + \tau \right), \tau_0 + k T^* + \tau \right) \cdot \prescript{*}{}{\bm{r}'} \left( s + \tau \right) \, \text{d}\tau} \right]}
        \nonumber \\
        =& \sum_{k = 0}^{a - 1}{\left[ \int_{0}^{T^*/2}{\bm{h}_2 \left( \prescript{*}{}{\bm{X}} \left( s + \tau \right), \tau_0 + k T^* + \tau \right) \cdot \prescript{*}{}{\bm{r}'} \left( s + \tau \right) \, \text{d}\tau} \right]}
        \label{eqn:HP2M3} \\
        &+ \sum_{k = 0}^{a - 1}{\left[ \int_{0}^{T^*/2}{\bm{h}_2 \left( \prescript{*}{}{\bm{X}} \left( s + T^* - \tau \right), \tau_0 + \left( k + 1 \right) T^* - \tau \right) \cdot \prescript{*}{}{\bm{r}'} \left( s + T^* - \tau \right) \, \text{d}\tau} \right]}
        \nonumber \\
        =& \sum_{k = 0}^{a - 1}{\left[
        \int_{0}^{T^*/2}{J \left( \prescript{*}{}{\bm{X}} \left( s + \tau \right) \right) \cos{\left( 2 \left( \tau_0 + k T^* + \tau \right) \right)} \, \text{d}\tau} \right.}
        \label{eqn:HP2M4} \\
        &\left. + \int_{0}^{T^*/2}{K \left( \prescript{*}{}{\bm{X}} \left( s + \tau \right) \right) \sin{\left( 2 \left( \tau_0 + k T^* + \tau \right) \right)} \, \text{d}\tau} \right.
        \nonumber \\
        &\left. - \int_{0}^{T^*/2}{\prescript{*}{}{z} \left( s + \tau \right) \prescript{*}{}{z'} \left( s + \tau \right) \, \text{d}\tau}
        \right]
        \nonumber \\
        &+ \sum_{k = 0}^{a - 1}{\left[
        \int_{0}^{T^*/2}{J \left( \prescript{*}{}{\bm{X}} \left( s + T^* - \tau \right) \right) \cos{\left( 2 \left( \tau_0 + \left( k + 1 \right) T^* - \tau \right) \right)} \, \text{d}\tau} \right.}
        \nonumber \\
        &\left. + \int_{0}^{T^*/2}{K \left( \prescript{*}{}{\bm{X}} \left( s + T^* - \tau \right) \right) \sin{\left( 2 \left( \tau_0 + \left( k + 1 \right) T^* - \tau \right) \right)} \, \text{d}\tau} \right.
        \nonumber \\
        &\left. - \int_{0}^{T^*/2}{\prescript{*}{}{z} \left( s + T^* - \tau \right) \prescript{*}{}{z'} \left( s + T^* - \tau \right) \, \text{d}\tau}
        \right]
        \nonumber
    \end{align}
    \label{eqn:HP2M}
\end{subequations}

\begin{subequations}
    \begin{align}
        \mathcal{M} \left( s, \tau_0 \right) =& \sum_{k = 0}^{a - 1}{\left[ C_k \right]} + \sum_{k = 0}^{a - 1}{\left[ D_k \right]}
        \label{eqn:HP2M5M} \\
        C_k =& \int_{0}^{T^*/2}{J \left( \prescript{*}{}{\bm{X}} \left( s + \tau \right) \right) \left( \cos{\left( 2 \left( \tau_0 + k T^* + \tau \right) \right)} - \cos{\left( 2 \left( \tau_0 + \left( k + 1 \right) T^* - \tau \right) \right)} \right) \, \text{d}\tau}
        \label{eqn:HP2M5Ck} \\
        D_k =& \int_{0}^{T^*/2}{K \left( \prescript{*}{}{\bm{X}} \left( s + \tau \right) \right) \left( \sin{\left( 2 \left( \tau_0 + k T^* + \tau \right) \right)} + \sin{\left( 2 \left( \tau_0 + \left( k + 1 \right) T^* - \tau \right) \right)} \right) \, \text{d}\tau}
        \label{eqn:HP2M5Dk}
    \end{align}
    \label{eqn:HP2M5}
\end{subequations}
We seek to evaluate both of the terms in \cref{eqn:HP2M5}, though we first develop expressions for some intermediate terms. Recall $a T^*$ must be an integer multiple of $\pi$ as $T_g = \pi$. Note that from trigonometry identities:
\begin{subequations}
    \begin{align}
        \sin{\left( 2 \left( \tau_0 + k T^* + \tau \right) \right)} + \sin{\left( 2 \left( \tau_0 + k T^* - \tau \right) \right)} &= 2 \sin{\left( 2 \left( \tau_0 + k T^* \right) \right)} \cos{\left( 2 \tau \right)}
        \label{eqn:trigid2s1} \\
        \cos{\left( 2 \left( \tau_0 + k T^* + \tau \right) \right)} - \cos{\left( 2 \left( \tau_0 + k T^* - \tau \right) \right)} &= -2 \sin{\left( 2 \left( \tau_0 + k T^* \right) \right)} \sin{\left( 2 \tau \right)}
        \label{eqn:trigid2c1} \\
        \sin{\left( 2 \left( \tau_0 + \tau \right) \right)} + \sin{\left( 2 \left( \tau_0 + a T^* - \tau \right) \right)} &= 2 \sin{ \left( 2 \tau_0 + a T^*\right)} \cos{\left( -a T^* + 2 \tau \right)}
        \label{eqn:trigid2s2} \\
        &= 2 \sin{ \left( 2 \tau_0 \right)} \cos{\left(2 \tau \right)}
        \nonumber \\
        \cos{\left( 2 \left( \tau_0 + \tau \right) \right)} - \cos{\left( 2 \left( \tau_0 + a T^* - \tau \right) \right)} &= -2 \sin{\left( 2 \tau_0 + a T^* \right)} \sin{\left( -a T^* + 2 \tau \right)}
        \label{eqn:trigid2c2} \\
        &= -2 \sin{ \left( 2 \tau_0 \right)} \sin{\left( 2 \tau \right)}
        \nonumber
    \end{align}
    \label{eqn:trigid2}
\end{subequations}
With psychic foresight we will present the following two equations. \Cref{eqn:HP2Jc,eqn:HP2Ks} are related to the terms in $C_k$ and $D_k$, respectively.
\begin{subequations}
    \begin{align}
        &\sum_{k = 0}^{a - 1}{\left[ \cos{\left( 2 \left( \tau_0 + k T^* + \tau \right) \right)} - \cos{\left( 2 \left( \tau_0 + \left( k + 1 \right) T^* - \tau \right) \right)} \right]} = \cos{\left( 2 \left( \tau_0 + 0 + \tau \right) \right)}
        \label{eqn:HP2Jc1} \\
        &\qquad \quad + \sum_{k = 1}^{a - 1}{\left[ \cos{\left( 2 \left( \tau_0 + k T^* + \tau \right) \right)} - \cos{\left( 2 \left( \tau_0 + k T^* - \tau \right) \right)} \right]} 
        \nonumber \\
        &\qquad \quad - \cos{\left( 2 \left( \tau_0 + a T^* - \tau \right) \right)}
        \nonumber \\
        &\qquad = -2 \sin{ \left( 2 \tau_0 \right)} \sin{\left( 2 \tau \right)} - \sum_{k = 1}^{a - 1}{\left[ 2 \sin{\left( 2 \left( \tau_0 + k T^* \right) \right)} \sin{\left( 2 \tau \right)} \right]}
        \label{eqn:HP2Jc2} \\
        &\qquad = -2 \sin{\left( 2 \tau \right)} \sum_{k = 0}^{a - 1}{\left[ \sin{\left( 2 \left( \tau_0 + k T^* \right) \right)} \right]}
        \label{eqn:HP2Jc3}
    \end{align}
    \label{eqn:HP2Jc}
\end{subequations}
\begin{subequations}
    \begin{align}
        &\sum_{k = 0}^{a - 1}{\left[ \sin{\left( 2 \left( \tau_0 + k T^* + \tau \right) \right)} + \sin{\left( 2 \left( \tau_0 + \left( k + 1 \right) T^* - \tau \right) \right)} \right]} = \sin{\left( 2 \left( \tau_0 + 0 + \tau \right) \right)}
        \label{eqn:HP2Ks1} \\
        &\qquad \quad + \sum_{k = 1}^{a - 1}{\left[ \sin{\left( 2 \left( \tau_0 + k T^* + \tau \right) \right)} + \sin{\left( 2 \left( \tau_0 + k T^* - \tau \right) \right)} \right]} 
        \nonumber \\
        &\qquad \quad + \sin{\left( 2 \left( \tau_0 + a T^* - \tau \right) \right)}
        \nonumber \\
        &\qquad = 2 \sin{ \left( 2 \tau_0 \right)} \cos{\left(2 \tau \right)} + \sum_{k = 1}^{a - 1}{\left[ 2 \sin{\left( 2 \left( \tau_0 + k T^* \right) \right)} \cos{\left( 2 \tau \right)} \right]}
        \label{eqn:HP2Ks2} \\
        &\qquad = 2 \cos{\left( 2 \tau \right)} \sum_{k = 0}^{a - 1}{\left[ \sin{\left( 2 \left( \tau_0 + k T^* \right) \right)} \right]}
        \label{eqn:HP2Ks3}
    \end{align}
    \label{eqn:HP2Ks}
\end{subequations}
Substituting the result in \cref{eqn:HP2Jc} back into \cref{eqn:HP2M5Ck} yields:
\begin{subequations}
    \begin{align}
        \sum_{k = 0}^{a - 1}{\left[ C_k \right]} &=
        \int_{0}^{T^*/2}{J \left( \prescript{*}{}{\bm{X}} \left( s + \tau \right) \right) \left( -2 \sin{\left( 2 \tau \right)} \sum_{k = 0}^{a - 1}{\left[ \sin{\left( 2 \left( \tau_0 + k T^* \right) \right)} \right]} \right) \, \text{d}\tau}
        \label{eqn:HP2J3} \\
        &=
        2 \left( \sum_{k = 0}^{a - 1}{\left[ \sin{\left( 2 \left( \tau_0 + k T^* \right) \right)} \right]} \right) \int_{0}^{T^*/2}{-J \left( \prescript{*}{}{\bm{X}} \left( s + \tau \right) \right) \sin{\left( 2 \tau \right)} \, \text{d}\tau}
        \label{eqn:HP2J4}
    \end{align}
    \label{eqn:HP2J}
\end{subequations}
Substituting the result in \cref{eqn:HP2Ks} back into \cref{eqn:HP2M5Dk} yields:
\begin{subequations}
    \begin{align}
        \sum_{k = 0}^{a - 1}{\left[ D_k \right]} &=
        \int_{0}^{T^*/2}{K \left( \prescript{*}{}{\bm{X}} \left( s + \tau \right) \right) \left( 2 \cos{\left( 2 \tau \right)} \sum_{k = 0}^{a - 1}{\left[ \sin{\left( 2 \left( \tau_0 + k T^* \right) \right)} \right]} \right) \, \text{d}\tau}
        \label{eqn:HP2K3} \\
        &=
        2 \left( \sum_{k = 0}^{a - 1}{\left[ \sin{\left( 2 \left( \tau_0 + k T^* \right) \right)} \right]} \right) \int_{0}^{T^*/2}{K \left( \prescript{*}{}{\bm{X}} \left( s + \tau \right) \right) \cos{\left( 2 \tau \right)} \, \text{d}\tau}
        \label{eqn:HP2K4}
    \end{align}
    \label{eqn:HP2K}
\end{subequations}

The result in \cref{eqn:HP2MSapp} can be obtained by substituting the expressions from \cref{eqn:HP2J} and \cref{eqn:HP2K} back into \cref{eqn:HP2M5} while using the following identity.
\begin{subequations}
    \begin{align}
        \sum_{k = 0}^{a - 1}{\left[ \sin{\left( 2 \left( \tau_0 + k T^* \right) \right)} \right]} &= 
        \csc{\left( T^* \right)} \sin{\left( a T^* \right)} \sin{\left( (a - 1) T^* + 2 \tau_0 \right)}
        \label{eqn:HP2ATIA} \\
        &= 
        \begin{cases}
          \sin{\left( 2 \tau_0 \right)}, & \text{if}\ a = 1 \\
          0, & \text{otherwise}
        \end{cases}
        \label{eqn:HP2ATIB}
    \end{align}
    \label{eqn:HP2ATI}
\end{subequations}

\end{proof}

If $\bm{h}_3$ had been used instead of $\bm{h_2}$ the same results shown in \cref{pro:P2} and \cref{pro:P3} would be obtained except:
\begin{subequations}
    \begin{align}
        \bm{h}_3 \left( \bm{X}, \alpha \right) \cdot \bm{r}' &= 
        J \left( \bm{X} \right) \cos{(2 \alpha)} + K \left( \bm{X} \right) \sin{(2 \alpha)} + 0
        \label{eqn:HP2h3rd1} \\
        J \left( \bm{X} \right) &= - \frac{38}{12} \left( 1 - \mu \right) \mu \left( -a x' + \left( y y' + z z' \right) \left( d \left( x + \mu \right) + e \right) \right)
        \label{eqn:HP2h3rdJ} \\
        K \left( \bm{X} \right) &= \frac{59}{12} \left( 1 - \mu \right) \mu \left( y x' \left( d \left( x + \mu \right) + e \right) - b y' + d y z z' \right)
        \label{eqn:HP2h3rdK}
    \end{align}
    \label{eqn:HP2h3rd}
\end{subequations}
While there are likely additional results relating to the Melnikov function that could be identified, however, the main results relevant to the analysis in this work--continuing CR3BP periodic orbits into the HR4BP--are presented and proved in the three Propositions of this appendix.


\clearpage

\bibliography{references}


\begin{thebibliography}{38}
\ifx \bisbn   \undefined \def \bisbn  #1{ISBN #1}\fi
\ifx \binits  \undefined \def \binits#1{#1}\fi
\ifx \bauthor  \undefined \def \bauthor#1{#1}\fi
\ifx \batitle  \undefined \def \batitle#1{#1}\fi
\ifx \bjtitle  \undefined \def \bjtitle#1{#1}\fi
\ifx \bvolume  \undefined \def \bvolume#1{\textbf{#1}}\fi
\ifx \byear  \undefined \def \byear#1{#1}\fi
\ifx \bissue  \undefined \def \bissue#1{#1}\fi
\ifx \bfpage  \undefined \def \bfpage#1{#1}\fi
\ifx \blpage  \undefined \def \blpage #1{#1}\fi
\ifx \burl  \undefined \def \burl#1{\textsf{#1}}\fi
\ifx \doiurl  \undefined \def \doiurl#1{\url{https://doi.org/#1}}\fi
\ifx \betal  \undefined \def \betal{\textit{et al.}}\fi
\ifx \binstitute  \undefined \def \binstitute#1{#1}\fi
\ifx \binstitutionaled  \undefined \def \binstitutionaled#1{#1}\fi
\ifx \bctitle  \undefined \def \bctitle#1{#1}\fi
\ifx \beditor  \undefined \def \beditor#1{#1}\fi
\ifx \bpublisher  \undefined \def \bpublisher#1{#1}\fi
\ifx \bbtitle  \undefined \def \bbtitle#1{#1}\fi
\ifx \bedition  \undefined \def \bedition#1{#1}\fi
\ifx \bseriesno  \undefined \def \bseriesno#1{#1}\fi
\ifx \blocation  \undefined \def \blocation#1{#1}\fi
\ifx \bsertitle  \undefined \def \bsertitle#1{#1}\fi
\ifx \bsnm \undefined \def \bsnm#1{#1}\fi
\ifx \bsuffix \undefined \def \bsuffix#1{#1}\fi
\ifx \bparticle \undefined \def \bparticle#1{#1}\fi
\ifx \barticle \undefined \def \barticle#1{#1}\fi
\bibcommenthead
\ifx \bconfdate \undefined \def \bconfdate #1{#1}\fi
\ifx \botherref \undefined \def \botherref #1{#1}\fi
\ifx \url \undefined \def \url#1{\textsf{#1}}\fi
\ifx \bchapter \undefined \def \bchapter#1{#1}\fi
\ifx \bbook \undefined \def \bbook#1{#1}\fi
\ifx \bcomment \undefined \def \bcomment#1{#1}\fi
\ifx \oauthor \undefined \def \oauthor#1{#1}\fi
\ifx \citeauthoryear \undefined \def \citeauthoryear#1{#1}\fi
\ifx \endbibitem  \undefined \def \endbibitem {}\fi
\ifx \bconflocation  \undefined \def \bconflocation#1{#1}\fi
\ifx \arxivurl  \undefined \def \arxivurl#1{\textsf{#1}}\fi
\csname PreBibitemsHook\endcsname

\bibitem[\protect\citeauthoryear{{\relax NASA Office of the Chief Financial Officer}}{2022}]{cit:NASASP22}
\begin{botherref}
\oauthor{\bsnm{{\relax NASA Office of the Chief Financial Officer}}}:
{NASA Strategic Plan 2022}.
\unskip\space N PD 1001 0D,
NASA
(Mar 2022)
\end{botherref}
\endbibitem

\bibitem[\protect\citeauthoryear{Whitley et~al.}{2018}]{cit:WhitleyDavisBurkeEtal}
\begin{bchapter}
\bauthor{\bsnm{Whitley}, \binits{R.J.}},
\bauthor{\bsnm{Davis}, \binits{D.C.}},
\bauthor{\bsnm{Burke}, \binits{L.M.}},
\bauthor{\bsnm{McCarthy}, \binits{B.P.}},
\bauthor{\bsnm{Power}, \binits{R.J.}},
\bauthor{\bsnm{Mcguire}, \binits{M.L.}},
\bauthor{\bsnm{Howell}, \binits{K.C.}}:
\bctitle{{Earth-Moon Near Rectilinear Halo and Butterfly Orbits for Lunar Surface Exploration}}.
In: \bbtitle{AAS/AIAA Astrodynamics Conference}
(\byear{2018})
\end{bchapter}
\endbibitem

\bibitem[\protect\citeauthoryear{Szebehely}{1967}]{cit:Szebehely}
\begin{bbook}
\bauthor{\bsnm{Szebehely}, \binits{V.}}:
\bbtitle{Theory of Orbits: The Restricted Problem of Three Bodies}.
\bpublisher{Academic Press Inc}\blocation{}
(\byear{1967})
\end{bbook}
\endbibitem

\bibitem[\protect\citeauthoryear{Park and Howell}{2022}]{cit:ParkHowell}
\begin{bchapter}
\bauthor{\bsnm{Park}, \binits{B.}},
\bauthor{\bsnm{Howell}, \binits{K.C.}}:
\bctitle{{Leveraging Intermediate Dynamical Models for Transitioning from the Circular Restricted Three-Body Problem to an Ephemeris Model}}.
In: \bbtitle{AAS/AIAA Astrodynamics Specialist Conference}
(\byear{2022})
\end{bchapter}
\endbibitem

\bibitem[\protect\citeauthoryear{Peng and Bai}{2018}]{cit:PengBai}
\begin{barticle}
\bauthor{\bsnm{Peng}, \binits{H.}},
\bauthor{\bsnm{Bai}, \binits{X.}}:
\batitle{{Natural deep space satellite constellation in the Earth-Moon elliptic system}}.
\bjtitle{Acta Astronautica}
\bvolume{153},
\bfpage{240}--\blpage{258}
(\byear{2018})
\end{barticle}
\endbibitem

\bibitem[\protect\citeauthoryear{Huang}{1960}]{cit:Huang}
\begin{botherref}
\oauthor{\bsnm{Huang}, \binits{S.-S.}}:
{Very Restricted Four-Body Problem}.
\unskip\space TN D-501,
{NASA}
(Sep 1960)
\end{botherref}
\endbibitem

\bibitem[\protect\citeauthoryear{G{\'o}mez et~al.}{2001}]{cit:GomezJorbaMasdemontEtal}
\begin{bchapter}
\bauthor{\bsnm{G{\'o}mez}, \binits{G.}},
\bauthor{\bsnm{Jorba}, \binits{{\'A}.}},
\bauthor{\bsnm{Masdemont}, \binits{J.}},
\bauthor{\bsnm{Sim{\'o}}, \binits{C.}}:
\bctitle{Normal form of the bicircular model and related topics}.
In: \bbtitle{Dynamics and Mission Design Near Libration Points},
pp. \bfpage{53}--\blpage{110}.
\bpublisher{World Scientific}\blocation{}
(\byear{2001})
\end{bchapter}
\endbibitem

\bibitem[\protect\citeauthoryear{Rosales}{2020}]{cit:Rosales}
\begin{botherref}
\oauthor{\bsnm{Rosales}, \binits{J.J.}}:
{On the effect of the Sun’s gravity around the Earth-Moon L1 and L2 libration points}.
PhD thesis,
Universitat de Barcelona
(2020)
\end{botherref}
\endbibitem

\bibitem[\protect\citeauthoryear{Jorba-Cusc{\'o} et~al.}{2018}]{cit:JorbaCuscoFarresJorba}
\begin{botherref}
\oauthor{\bsnm{Jorba-Cusc{\'o}}, \binits{M.}},
\oauthor{\bsnm{Farr{\'e}s}, \binits{A.}},
\oauthor{\bsnm{Jorba}, \binits{{\'A}.}}:
{Two Periodic Models for the Earth-Moon System}.
Frontiers in Applied Mathematics and Statistics
\textbf{4}
(2018)
\end{botherref}
\endbibitem

\bibitem[\protect\citeauthoryear{Sim{\'o} et~al.}{1995}]{cit:SimoGomezJorbaEtal}
\begin{bchapter}
\bauthor{\bsnm{Sim{\'o}}, \binits{C.}},
\bauthor{\bsnm{G{\'o}mez}, \binits{G.}},
\bauthor{\bsnm{Jorba}, \binits{{\'A}.}},
\bauthor{\bsnm{Masdemont}, \binits{J.}}:
\bctitle{{The Bicircular Model Near the Triangular Libration Points of the RTBP}}.
In: \bbtitle{From Newton to Chaos: Modern Techniques for Understanding and Coping with Chaos in N-Body Dynamical Systems},
pp. \bfpage{343}--\blpage{370}.
\bpublisher{Springer}\blocation{}
(\byear{1995})
\end{bchapter}
\endbibitem

\bibitem[\protect\citeauthoryear{Boudad et~al.}{2020}]{cit:BoudadHowellDavis}
\begin{barticle}
\bauthor{\bsnm{Boudad}, \binits{K.K.}},
\bauthor{\bsnm{Howell}, \binits{K.C.}},
\bauthor{\bsnm{Davis}, \binits{D.C.}}:
\batitle{{Dynamics of synodic resonant near rectilinear halo orbits in the bicircular four-body problem}}.
\bjtitle{Advances in Space Research}
\bvolume{66}(\bissue{9}),
\bfpage{2194}--\blpage{2214}
(\byear{2020})
\end{barticle}
\endbibitem

\bibitem[\protect\citeauthoryear{Castell{\'a} and Jorba}{2000}]{cit:CastellaJorba}
\begin{barticle}
\bauthor{\bsnm{Castell{\'a}}, \binits{E.}},
\bauthor{\bsnm{Jorba}, \binits{{\'A}.}}:
\batitle{{On the vertical families of two-dimensional tori near the triangular points of the Bicircular problem}}.
\bjtitle{Celestial Mechanics and Dynamical Astronomy}
\bvolume{76},
\bfpage{35}--\blpage{54}
(\byear{2000})
\end{barticle}
\endbibitem

\bibitem[\protect\citeauthoryear{Rosales et~al.}{2021}]{cit:RosalesJorbaJorbaCuscoBCP}
\begin{botherref}
\oauthor{\bsnm{Rosales}, \binits{J.J.}},
\oauthor{\bsnm{Jorba}, \binits{{\'A}.}},
\oauthor{\bsnm{Jorba-Cusc{\'o}}, \binits{M.}}:
{Families of Halo-like invariant tori around L2 in the Earth-Moon Bicircular Problem}.
Celestial Mechanics and Dynamical Astronomy
\textbf{133}
(2021)
\end{botherref}
\endbibitem

\bibitem[\protect\citeauthoryear{Andreu}{1998}]{cit:Andreu}
\begin{botherref}
\oauthor{\bsnm{Andreu}, \binits{M.A.}}:
{The Quasi-bicircular Problem}.
PhD thesis,
Universitat de Barcelona
(1998)
\end{botherref}
\endbibitem

\bibitem[\protect\citeauthoryear{Rosales et~al.}{2023}]{cit:RosalesJorbaJorbaCuscoQBCP}
\begin{botherref}
\oauthor{\bsnm{Rosales}, \binits{J.J.}},
\oauthor{\bsnm{Jorba}, \binits{{\'A}.}},
\oauthor{\bsnm{Jorba-Cusc{\'o}}, \binits{M.}}:
{Invariant manifolds near L1 and L2 in the quasi-bicircular problem}.
Celestial Mechanics and Dynamical Astronomy
\textbf{135}
(2023)
\end{botherref}
\endbibitem

\bibitem[\protect\citeauthoryear{Scheeres}{1998}]{cit:ScheeresHR4BP}
\begin{barticle}
\bauthor{\bsnm{Scheeres}, \binits{D.J.}}:
\batitle{{The Restricted Hill Four-Body Problem with Applications to the Earth–Moon–Sun System}}.
\bjtitle{Celestial Mechanics and Dynamical Astronomy}
\bvolume{70}(\bissue{2}),
\bfpage{75}--\blpage{98}
(\byear{1998})
\end{barticle}
\endbibitem

\bibitem[\protect\citeauthoryear{Peterson et~al.}{2023}]{cit:PetersonRosalesScheeres}
\begin{botherref}
\oauthor{\bsnm{Peterson}, \binits{L.T.}},
\oauthor{\bsnm{Rosales}, \binits{J.J.}},
\oauthor{\bsnm{Scheeres}, \binits{D.J.}}:
{The vicinity of Earth–Moon L1 and L2 in the Hill restricted 4-body problem}.
Physica D: Nonlinear Phenomena
\textbf{455}
(2023)
\end{botherref}
\endbibitem

\bibitem[\protect\citeauthoryear{Olikara et~al.}{2016}]{cit:OlikaraGomezMasdemont}
\begin{bchapter}
\bauthor{\bsnm{Olikara}, \binits{Z.P.}},
\bauthor{\bsnm{G{\'o}mez}, \binits{G.}},
\bauthor{\bsnm{Masdemont}, \binits{J.J.}}:
\bctitle{{A Note on Dynamics About the Coherent Sun-Earth-Moon Collinear Libration Points}}.
In: \beditor{\bsnm{G{\'o}mez}, \binits{G.}},
\beditor{\bsnm{Masdemont}, \binits{J.J.}} (eds.)
\bbtitle{Astrodynamics Network AstroNet-II},
pp. \bfpage{183}--\blpage{192}.
\bpublisher{Springer}\blocation{}
(\byear{2016})
\end{bchapter}
\endbibitem

\bibitem[\protect\citeauthoryear{Henry et~al.}{2023a}]{cit:HenryRosalesBrownEtalISTS}
\begin{bchapter}
\bauthor{\bsnm{Henry}, \binits{D.}},
\bauthor{\bsnm{Rosales}, \binits{J.}},
\bauthor{\bsnm{Brown}, \binits{G.}},
\bauthor{\bsnm{Peterson}, \binits{L.}},
\bauthor{\bsnm{Scheeres}, \binits{D.}}:
\bctitle{{Quasi-Periodic Orbits near Earth-Moon L1 in the Hill Restricted Four-Body Problem}}.
In: \bbtitle{34th ISTS}
(\byear{2023})
\end{bchapter}
\endbibitem

\bibitem[\protect\citeauthoryear{Henry et~al.}{2023b}]{cit:HenryRosalesBrownEtalAASC}
\begin{bchapter}
\bauthor{\bsnm{Henry}, \binits{D.B.}},
\bauthor{\bsnm{Rosales}, \binits{J.}},
\bauthor{\bsnm{Brown}, \binits{G.M.}},
\bauthor{\bsnm{Scheeres}, \binits{D.J.}}:
\bctitle{{Quasi-Periodic Orbits near Earth-Moon L1 and L2 in the Hill Restricted Four-Body Problem}}.
In: \bbtitle{AAS/AIAA Astrodynamics Specialist Conference}
(\byear{2023})
\end{bchapter}
\endbibitem

\bibitem[\protect\citeauthoryear{Peterson et~al.}{2024}]{cit:PetersonJorbaBrownEtal}
\begin{botherref}
\oauthor{\bsnm{Peterson}, \binits{L.T.}},
\oauthor{\bsnm{Jorba}, \binits{A.}},
\oauthor{\bsnm{Brown}, \binits{G.M.}},
\oauthor{\bsnm{Scheeres}, \binits{D.J.}}:
{Dynamics Around the Earth-Moon Triangular Points in the Hill Restricted 4-Body Problem}.
Communications in Nonlinear Science and Numerical Simulation
(2024).
In Preparation
\end{botherref}
\endbibitem

\bibitem[\protect\citeauthoryear{Olikara and Scheeres}{2017}]{cit:OlikaraScheeres}
\begin{bchapter}
\bauthor{\bsnm{Olikara}, \binits{Z.P.}},
\bauthor{\bsnm{Scheeres}, \binits{D.J.}}:
\bctitle{{Mapping Connections Between Planar Sun-Earth-Moon Libration Orbits}}.
In: \bbtitle{27th AAS/AIAA Space Flight Mechanics Meeting}
(\byear{2017})
\end{bchapter}
\endbibitem

\bibitem[\protect\citeauthoryear{Sanaga and Howell}{2023}]{cit:SanagaHowell}
\begin{bchapter}
\bauthor{\bsnm{Sanaga}, \binits{R.R.}},
\bauthor{\bsnm{Howell}, \binits{K.C.}}:
\bctitle{{Synodic Resonant Halo Orbits in the Hill Restricted Four-Body Problem}}.
In: \bbtitle{33rd AAS/AIAA Space Flight Mechanics Meeting}
(\byear{2023})
\end{bchapter}
\endbibitem

\bibitem[\protect\citeauthoryear{Melnikov}{1963}]{cit:Melnikov}
\begin{barticle}
\bauthor{\bsnm{Melnikov}, \binits{V.K.}}:
\batitle{{On the stability of the center for time periodic perturbations}}.
\bjtitle{Transactions of the Moscow Mathematical Society}
\bvolume{12},
\bfpage{1}--\blpage{57}
(\byear{1963})
\end{barticle}
\endbibitem

\bibitem[\protect\citeauthoryear{Greenspan and Holmes}{1981}]{cit:GreenspanHolmes}
\begin{botherref}
\oauthor{\bsnm{Greenspan}, \binits{B.D.}},
\oauthor{\bsnm{Holmes}, \binits{P.}}:
{Homoclinic Orbits, Subharmonics and Global Bifurcations in Forced Oscillations}.
\unskip\space ADA103564,
{U.S. Army Research Office}
(Jun 1981)
\end{botherref}
\endbibitem

\bibitem[\protect\citeauthoryear{Wiggins}{2003}]{cit:Wiggins}
\begin{bbook}
\bauthor{\bsnm{Wiggins}, \binits{S.}}:
\bbtitle{Introduction to Applied Nonlinear Dynamical Systems and Chaos},
\bedition{2nd} edn.
\bsertitle{Texts in Applied Mathematics}.
\bpublisher{Springer}\blocation{}
(\byear{2003})
\end{bbook}
\endbibitem

\bibitem[\protect\citeauthoryear{Guckenheimer and Holmes}{1983}]{cit:GuckenheimerHolmes}
\begin{bbook}
\bauthor{\bsnm{Guckenheimer}, \binits{J.}},
\bauthor{\bsnm{Holmes}, \binits{P.}}:
\bbtitle{Nonlinear Dylanicas, Systems, and Bifurcations of Vector Fields}.
\bsertitle{Applied Mathematical Sciences}.
\bpublisher{Springer}\blocation{}
(\byear{1983})
\end{bbook}
\endbibitem

\bibitem[\protect\citeauthoryear{Perko}{2001}]{cit:Perko}
\begin{bbook}
\bauthor{\bsnm{Perko}, \binits{L.}}:
\bbtitle{Differential Equations and Dynamical Systems},
\bedition{3rd} edn.
\bsertitle{Texts in Applied Mathematics}.
\bpublisher{Springer}\blocation{}
(\byear{2001})
\end{bbook}
\endbibitem

\bibitem[\protect\citeauthoryear{Haller}{1999}]{cit:Haller}
\begin{bbook}
\bauthor{\bsnm{Haller}, \binits{G.}}:
\bbtitle{Chaos Near Resonance}.
\bsertitle{Applied Mathematical Sciences}.
\bpublisher{Springer}\blocation{}
(\byear{1999})
\end{bbook}
\endbibitem

\bibitem[\protect\citeauthoryear{Guo et~al.}{2022}]{cit:GuoTianXueEtal}
\begin{botherref}
\oauthor{\bsnm{Guo}, \binits{X.}},
\oauthor{\bsnm{Tian}, \binits{R.}},
\oauthor{\bsnm{Xue}, \binits{Q.}},
\oauthor{\bsnm{Zhang}, \binits{X.}}:
{Sub-harmonic Melnikov function for a high-dimensional non-smooth coupled system}.
Chaos, Solitons \& Fractals
\textbf{164}
(2022)
\end{botherref}
\endbibitem

\bibitem[\protect\citeauthoryear{Veerman and Holmes}{1985}]{cit:VeermanHolmes}
\begin{barticle}
\bauthor{\bsnm{Veerman}, \binits{P.}},
\bauthor{\bsnm{Holmes}, \binits{P.}}:
\batitle{The existence of arbitrarily many distinct periodic orbits in a two degree of freedom hamiltonian system}.
\bjtitle{Physica D: Nonlinear Phenomena}
\bvolume{14}(\bissue{2}),
\bfpage{177}--\blpage{192}
(\byear{1985})
\end{barticle}
\endbibitem

\bibitem[\protect\citeauthoryear{Yagasaki}{1996}]{cit:Yagasaki}
\begin{barticle}
\bauthor{\bsnm{Yagasaki}, \binits{K.}}:
\batitle{The melnikov theory for subharmonics and their bifurcations in forced oscillations}.
\bjtitle{SIAM Journal on Applied Mathematics}
\bvolume{56}(\bissue{6}),
\bfpage{1720}--\blpage{1765}
(\byear{1996})
\end{barticle}
\endbibitem

\bibitem[\protect\citeauthoryear{Polcar and Semer{\'a}k}{2019}]{cit:PolcarSemerak}
\begin{botherref}
\oauthor{\bsnm{Polcar}, \binits{L.}},
\oauthor{\bsnm{Semer{\'a}k}, \binits{O.}}:
{Free motion around black holes with discs or rings: Between integrability and chaos. VI. The Melnikov method}.
Phys. Rev. D
\textbf{100}(10)
(2019)
\end{botherref}
\endbibitem

\bibitem[\protect\citeauthoryear{Cenedese and Haller}{2020}]{cit:CenedeseHaller}
\begin{botherref}
\oauthor{\bsnm{Cenedese}, \binits{M.}},
\oauthor{\bsnm{Haller}, \binits{G.}}:
{How do conservative backbone curves perturb into forced responses? A Melnikov function analysis}.
Proceedings of the Royal Society A
\textbf{476}(2234)
(2020)
\end{botherref}
\endbibitem

\bibitem[\protect\citeauthoryear{Rhouma and Chicone}{2000}]{cit:RhoumaChicone}
\begin{barticle}
\bauthor{\bsnm{Rhouma}, \binits{M.B.H.}},
\bauthor{\bsnm{Chicone}, \binits{C.}}:
\batitle{{On the Continuation of Periodic Orbits}}.
\bjtitle{Methods and Applications of Analysis}
\bvolume{7}(\bissue{1}),
\bfpage{85}--\blpage{104}
(\byear{2000})
\end{barticle}
\endbibitem

\bibitem[\protect\citeauthoryear{Scheeres}{2012}]{cit:ScheeresBook}
\begin{bbook}
\bauthor{\bsnm{Scheeres}, \binits{D.J.}}:
\bbtitle{Orbital Motion in Strongly Perturbed Environments}.
\bpublisher{Springer}\blocation{}
(\byear{2012})
\end{bbook}
\endbibitem

\bibitem[\protect\citeauthoryear{Holmes}{1990}]{cit:Holmes}
\begin{barticle}
\bauthor{\bsnm{Holmes}, \binits{P.}}:
\batitle{{Poincar{\'e}, celestial mechanics, dynamical-systems theory and ``chaos''}}.
\bjtitle{Physics Reports}
\bvolume{193}(\bissue{3}),
\bfpage{137}--\blpage{163}
(\byear{1990})
\end{barticle}
\endbibitem

\bibitem[\protect\citeauthoryear{Wintner}{1952}]{cit:Wintner}
\begin{bbook}
\bauthor{\bsnm{Wintner}, \binits{A.}}:
\bbtitle{The Analytical Foundations of Celestial Mechanics}.
\bsertitle{Dover Books on Physics}.
\bpublisher{Dover Publications}\blocation{}
(\byear{1952})
\end{bbook}
\endbibitem

\end{thebibliography}

\end{document}